%------------------------------------------------------------------------------
% nakayama2.tex 
%                    
% Beginn: 01.07.2012, corrected: 30.07.2012, Neufassung nach Referentenbericht
%------------------------------------------------------------------------------
\magnification=\magstep1   
\input amstex
\UseAMSsymbols
\input pictex
%\hoffset=0truecm \voffset=0truecm 
\vsize=23truecm
\NoBlackBoxes
\parindent=18pt
  
   \font\rmk=cmr8    \font\itk=cmti8  \font\ttk=cmtt8

%\hrule height 2pt \vskip 3pt \hrule \bigskip\bigskip\bigskip

\def\mo{\operatorname{mod}}
\def\gp{\operatorname{gp}}

\def\Hom{\operatorname{Hom}}

\def\Ext{\operatorname{Ext}}

\def\rad{\operatorname{rad}}
\def\add{\operatorname{add}}
\def\Ker{\operatorname{Ker}}

\def\CM{\operatorname{CM}}
\def\soc{\operatorname{soc}}
\def\top{\operatorname{top}}
\def\Tr{\operatorname{Tr}}
\def\Im{\operatorname{Im}}

\def\arr#1#2{\arrow <1.5mm> [0.25,0.75] from #1 to #2}

%%%%%%%%%%%%%%%%%%%%%%%%%%%%%%%%%%%%%%%%%%%%%%%%%%%%%%%%%%%%%%
	\vglue1truecm
\plainfootnote{}
{\rmk 2010 \itk Mathematics Subject Classification. \rmk 
Primary 
        16G10, %Representations of artinian rings 
        16G50. % Cohen-Macaulay modules.
Secondary:
        16D90, %Module categories [See also 16Gxx, 16S90]; module theory in a 
               % category-theoretic context; Morita equivalence and duality
        16E05, % Syzygies, resolutions, complexes 
        16G70, % Auslander-Reiten sequences (almost split sequences)
                  %and  Auslander-Reiten quivers
        18G25. % Relative homological algebra
}
\centerline{\bf The Gorenstein projective modules for the Nakayama algebras. I.}
                     \bigskip
\centerline{Claus Michael Ringel}     
		  	    \bigskip\medskip

{\narrower\narrower
Abstract: The aim of this paper is to outline the structure of the category of 
the Gorenstein projective
$\Lambda$-modules, where $\Lambda$ is a Nakayama algebra.  
In addition, we are going to introduce 
the resolution quiver of $\Lambda$. It provides a fast algorithm 
in order to obtain the Gorenstein projective $\Lambda$-modules 
and to decide
whether $\Lambda$ is a Gorenstein algebra or not, and whether it is $\CM$-free or not.
\par} 
      \bigskip
Throughout the paper, $\Lambda$ will be a  Nakayama algebra without simple projective modules, and 
the modules will be left $\Lambda$-modules of finite length.
Let $\mo\Lambda$ be the category of all such modules. 
The subcategories of $\mo\Lambda$ which we will deal with 
will be assumed to be full and
closed under direct sums and direct summands. If $\Cal M$ is a class of modules,
we denote by $\add\Cal M$ the smallest subcategory containing $\Cal M$
and $(\add\Cal M)$-approximations will be just called $\Cal M$-approximations. 

We denote by $\gp \Lambda$ the full subcategory of $\mo\Lambda$ given by
the Gorenstein projective modules, and $\gp_0\Lambda$ denotes the subcategory
of all Gorenstein projective modules without any indecomposable projective
direct summand. Recall that $\Lambda$ is said to be {\it $\CM$-free} [C] provided all
Gorenstein projective modules are projective, thus provided $\gp_0\Lambda$ is
the zero category. Our aim is to
describe the subcategory $\Cal C = \Cal C(\Lambda)$ whose indecomposable objects are the indecomposable 
non-projective Gorenstein projective modules as well as their projective covers.
We call $\Cal C$ the {\it Gorenstein core,} clearly
$$
 \gp_0\Lambda\ \subseteq\ \Cal C \ \subseteq\ \gp\Lambda.
$$
The first assertions concern the structure of $\Cal C(\Lambda).$ Here, 
$\underline {\gp}\ \Lambda$ denotes the factor category of $\gp \Lambda$ obtained
by factoring out the ideal of all maps which factor through a projective module. 
Similarly, $\underline {\mo}\ \Lambda'$ is obtained from $\mo \Lambda'$ by 
factoring out the ideal of all maps which factor through a projective 
$\Lambda'$-module. 
		   \bigskip
%%%%%%%%%%%%%%%%%%%%%%%%%%%%%%
{\bf Proposition 1.} {\it Let $\Lambda$ be a  Nakayama algebra. 
The Gorenstein core
$\Cal C = \Cal C(\Lambda)$ is a full exact abelian subcategory of $\mo \Lambda$ which is
closed under extensions, projective covers and minimal left $\Lambda$-approximations;
$\Cal C$ is equivalent to $\mo\Lambda',$ 
where $\Lambda'$ is a self-injective Nakayama algebra 
and the inclusion functor $\Cal C \to \gp \Lambda$ induces an
equivalence $\underline {\mo}\ \Lambda' \to \underline{\gp}\ \Lambda$.
If $\Cal C$ is not zero and $\Lambda$ is connected, then also $\Lambda'$ is connected.}
   \bigskip
Given a class $\Cal M$ of modules, we denote by $\Cal F(\Cal M)$ the class of 
all modules with a filtration with factors in $\Cal M$. 
    \medskip	 
%%%%%%%%%%%%%%%%%%%%%%%%%%%%%%%%%%%%%%%%%
{\bf Proposition 2.} {\it Let $\Lambda$ be a 
Nakayama algebra. Let $\Cal E = \Cal E(\Lambda)$ 
be the
class of non-zero modules in $\gp_0\Lambda$ such that no proper non-zero factor module is Gorenstein projective.

{\rm (a)} We have $\Cal C(\Lambda) = \Cal F(\Cal E)$.

{\rm (b)} If $\Lambda$ is connected and $\Cal C(\Lambda)$ is not zero, 
let $E_1,\dots, E_g$ be  representatives of the isomorphism classes in $\Cal E$.  
Then any simple $\Lambda$-module occurs with multiplicity $1$ in
$\bigoplus_{i=1}^g E_i$. In particular, $E_1,\dots, E_g$ are pairwise orthogonal bricks.

{\rm (c)}  
Let $\Cal E'$ be the
class of non-zero modules in $\gp_0\Lambda$ 
such that no proper non-zero submodule is Gorenstein projective.
Then $\Cal E = \Cal E'$.}
     \medskip
Since we deal with a set of orthogonal bricks, the elements of  
$\Cal E(\Lambda)$ are just the simple objects of $\Cal C(\Lambda)$ (see
for example [R1]), we call them the {\it elementary} Gorenstein projective modules.
Note that assertion (a) implies that {\it $\Cal C(\Lambda)$ is the extension closure
of $\gp_0$} (since $\Cal E \subseteq \gp_0 \subseteq \Cal C(\Lambda)).$  
   \bigskip 
Given a module $M$, we denote by $P(M)$ the projective cover, by $\Omega(M)$ the
first syzygy module  and by $\tau M = D\Tr M,\ 
\tau^-M = \Tr D M$ the Auslander-Reiten translates of $M$.
If $\Cal M$ is a class of modules, we write
$\tau \Cal M$ for the class of modules $\tau M$ with $M\in \Cal M$, and similarly,
$\tau^- \Cal M$ is the class of modules $\tau^- M$ with $M\in \Cal M$.
	\medskip
{\bf Proposition 3.}  {\it Let $\Lambda$ be a Nakayama algebra.
Let $\Cal X = \Cal X(\Lambda)$ be the class of simple modules $S$ with $P(S)$ belonging
to $\Cal C(\Lambda).$ The following conditions are equivalent for a module $M$:
\item{\rm(i)} $M$ belongs to $\Cal C(\Lambda).$
\item{\rm(ii)} $\top M$ belongs to $\add \Cal X$, and $\soc M$ belongs to $\add\tau^-\Cal X$.
\item{\rm(iii)} $\top M$ and $\top \Omega M$ both belong to $\add \Cal X.$\par}
		\medskip
We see that $\Cal C(\Lambda)$ may be obtained by deleting ray and corays from the Auslander-Reiten
quiver of $\Lambda$. Namely, we have to delete the rays consisting of the indecomposable 
modules with top not in $\Cal X$, as well as the corays  consisting of the indecomposable 
modules with socle not in $\tau^-\Cal X$.
	\bigskip 
%=====================================================
The basic observation which we use is the following characterization of the
indecomposable non-projective Gorenstein projective modules. Here, an indecomposable
projective module $P$ is said to be {\it minimal projective} provided no proper non-zero submodule
of $P$ is projective, or, equivalently, provided the projective dimension of $\top P$
is at least $2$.
   \medskip
{\bf Proposition 4.} {\it Let $\Lambda$ be a Nakayama algebra. Let 
$M$ be an indecomposable non-projective module.
The following assertions are equivalent:
\item{\rm (i)} The module $M$ is Gorenstein projective.
\item{\rm (ii)} All the projective modules occurring in a minimal projective resolution of
$M$ are minimal projective.
\item{\rm (iii)} There is an exact sequence
$$
 0 \to M \to P_{n-1} \to \cdots \to P_0 \to M \to 0 \tag{$*$}
$$
such that all the modules $P_i$ are minimal projective.\par}
     \bigskip
Of course, since any Nakayama algebra is representation-finite, any Gorenstein projective
module $M$ has a periodic projective resolution, thus there is an exact sequence 
of the form $(*)$
with projective modules $P_0,\dots P_{n-1}$. But usually not all modules with a periodic
projective resolution are Gorenstein projective. The interesting feature here is the fact
that the isomorphism classes of the modules $P_i$ in a periodic projective resolution of $M$
determine whether $M$ is Gorenstein projective or not.
	  \bigskip 
Throughout the paper, we will denote the number of simple modules by $s = s(\Lambda)$,
the minimal length of an indecomposable projective module by $p = p(\Lambda)$. Since we assume that
there are no simple projective modules, we have $p\ge 2.$

Note that $p > s$
if and only if no projective module is a brick. Examples which we will exhibit 
in section 7 show
that algebras with $p \le s$ may have some irregularities, thus some of our results
require the condition $p > s$. 

	\medskip
We will work with the {\it resolution quiver} $R = R(\Lambda)$ of $\Lambda$:
The vertices of  $R$ are the simple modules and for every vertex $S$, 
there is an arrow from $S$ to $\tau\soc P(S)$. 
Since any vertex in $R$ is the start of a unique arrow, 
any connected component of $R$
contains precisely one cycle. A vertex $S$ of $R$ is
said to be {\it black} provided the projective dimension of $S$ is at least $2$ (thus if and
only if $P(S)$ is a minimal projective module), otherwise it is said to be {\it red.}
As we will see 
the modules in $\Cal X(\Lambda)$ are precisely the simple modules
$S$ which belong to a cycle of black vertices provided $p > s,$ see Corollary 3 to Lemma 5.
    \medskip 
{\bf Proposition 5.} {\it Let $\Lambda$ be a Nakayama algebra and assume that $p > s$.

{\rm (a)}  The algebra $\Lambda$ is a Gorenstein algebra if and only if  
any cycle in $R(\Lambda)$ contains only black vertices. 

{\rm (b)}  The algebra $\Lambda$ is $\CM$-free if and only if any cycle 
in $R(\Lambda)$ contains at least one red vertex.}
   \bigskip
In a second part [R2] we will describe some further properties of the
resolution quiver. In particular, we will show that for a connected Nakayama algebra $\Lambda$ 
there is either no loop in $R(\Lambda)$, or else all cycles in $R(\Lambda)$ are loops.
This result has also been obtained (independently and with a different proof) by
Dawei Shen [S].
      \bigskip\bigskip
%==============================================================
{\bf An example.} Let $Q$ be the quiver of type $\widetilde{\Bbb A}_4$ with cyclic orientation, say with
vertices $1,2,3,4,5$ and arrows $i\to i\!+\!1$ (modulo $5$). Since a Nakayama algebra
is defined by zero relations (monomials), it is
sufficient to mention the length of the indecomposable projectives, 
instead of writing down the relations. 
Here we consider the case
where $|P(i)| = 13$ for $i=1,2$ and $|P(i)|= 12$ for $i=3,4,5$, thus we deal with the algebra
$\Lambda$ with Kupisch series $(13,13,12,12,12)$. There are two elementary
Gorenstein projective modules, namely $E(1)$ with composition factors $1,2,3$ and
$E(4)$ with composition factors $4,5$. The right picture below shows the
support of these modules $E(1)$ and $E(4)$, thus the $\Ext$-quiver for
$\Cal E = \{E(1),E(4)\}$ will be the quiver of type $\widetilde{\Bbb A}_1$ with
cyclic orientation.
$$
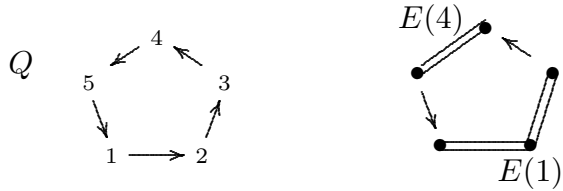
\hbox{\beginpicture
\setcoordinatesystem units <1cm,1cm>
%============================================
\put{\beginpicture
\put{$Q$} at -1 1
\put{$\ssize 1$} at 0 0 
\put{$\ssize 5$} at -.25 .8
\put{$\ssize 4$} at 0.5 1.3
\put{$\ssize 3$} at 1.25 .8
\put{$\ssize 2$} at 1 0
\arr{0.2 0}{0.8 0}
\arr{-.2 .6}{-.05 .2}
\arr{1.05 .2}{1.2 .6}
\arr{0.3 1.2}{0 1}
\arr{1 1}{0.7 1.2}

\endpicture} at 0 0
%============================================
\put{\beginpicture
%\put{$Q'$} at -1.2 1
\put{$\bullet$} at 0 0 
\put{$\bullet$} at -.25 .8
\put{$\bullet$} at 0.5 1.3
\put{$\bullet$} at 1.25 .8
\put{$\bullet$} at 1 0
\plot 0 0.05  1 0.05 /
\plot 0 -.04  1 -.04 /
\arr{-.2 .6}{-.05 .2}
\plot 1.05 0  1.3 .8 /
\plot 0.95 0  1.2 .8 /
\arr{1 1}{0.7 1.2}
\plot -.2 .8  0.5 1.3 /
\plot -.2 .9  0.5 1.4 /
\put{$E(1)$} at 1 -.3
\put{$E(4)$} at -.1 1.4

\endpicture} at 4 0
\endpicture}
$$

Next, let us draw twice the Auslander-Reiten
quiver of $\Lambda$
(always, the left dashed boundary has to be identified with the right dashed boundary)
and label the simple module $S(i)$ just by $i$. 
On the left,
we use bullets to mark the indecomposable objects in $\Cal C(\Lambda)$;
the encircled ones
are the elementary Gorenstein projective modules $E(1)$ and $E(4)$).
Then in the middle, we shade the rays and the
corays which have to be deleted: note that 
$\Cal X(\Lambda) = \{S(1),S(4)\},$ thus $\tau^-\Cal X(\Lambda) = \{S(5),S(3)\},$ this means that
we have to delete the corays ending in $S(2), S(3),S(5)$, and the rays starting in 
$S(1), S(2), S(4)$. On the right, we present the Auslander-Reiten quiver of $\Cal C(\Lambda)$.

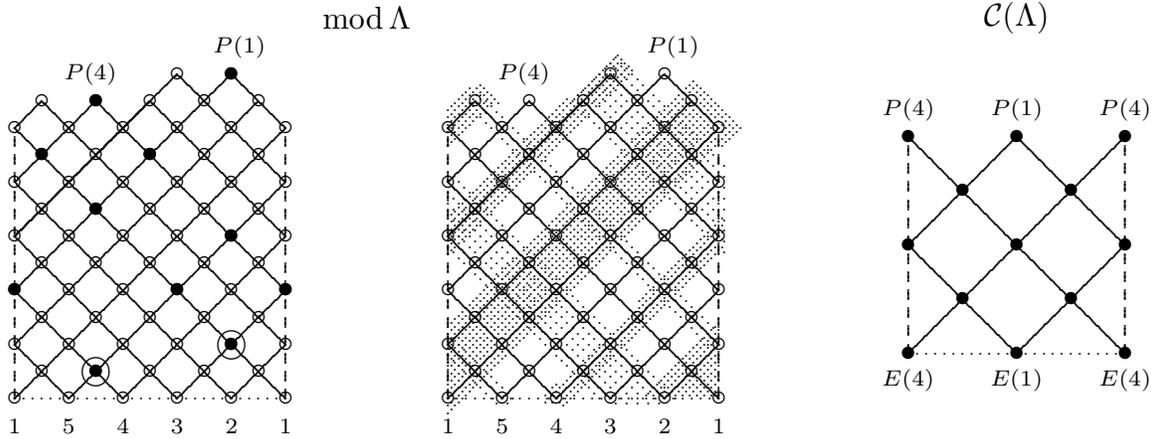
$$
\hbox{\beginpicture
\setcoordinatesystem units <.3cm,.3cm>
%============================================
\put{\beginpicture
\multiput{$\circ$} at 
   0 0  0 2  0 4  0 6  0 8  0 10 
   1 1  1 3  1 5  1 7  1 9  1 11
   2 0  2 2  2 4  2 6  2 8  2 10 
   3 1  3 3  3 5  3 7  3 9  3 11
   4 0  4 2  4 4  4 6  4 8  4 10 
   5 1  5 3  5 5  5 7  5 9  5 11
   6 0  6 2  6 4  6 6  6 8  6 10  6 12
   7 1  7 3  7 5  7 7  7 9  7 11
   8 0  8 2  8 4  8 6  8 8  8 10  8 12
   9 1  9 3  9 5  9 7  9 9  9 11
   10 0  10 2  10 4  10 6  10 8  10 10 /
\plot 0 0  10 10  8 12  0 4  4 0  10 6  5 11  
   0 6  6 12  10 8  2 0  0 2  9 11 /
\plot  8 0   10 2  1 11  0 10  10 0  /
\plot 0 6  6 0  10 4  3 11  0 8  8 0 /
\setdashes <1mm>
\plot 0 0  0 10 /
\plot 10 0  10 10 /
\setdots <1mm> 
\plot 0 0  10 0 /
\multiput{$\bullet$} at 0 4  1 9  3 1  3 7  3 11  5 9  6 4  8 2  8 6  8 12  10 4 /
\multiput{$\bigcirc$} at 3 1  8 2 /
\put{$\ssize P(1)$} at 8.3 13
\put{$\ssize P(4)$} at 2.8 12
\put{$\ssize 1$} at 10 -1
\put{$\ssize 2$} at 8 -1
\put{$\ssize 3$} at 6 -1
\put{$\ssize 4$} at 4 -1
\put{$\ssize 5$} at 2 -1
\put{$\ssize 1$} at 0 -1
\endpicture} at 0 0 
%============================================
\put{\beginpicture
\multiput{$\circ$} at 
   0 0  0 2  0 4  0 6  0 8  0 10 
   1 1  1 3  1 5  1 7  1 9  1 11
   2 0  2 2  2 4  2 6  2 8  2 10 
   3 1  3 3  3 5  3 7  3 9  3 11
   4 0  4 2  4 4  4 6  4 8  4 10 
   5 1  5 3  5 5  5 7  5 9  5 11
   6 0  6 2  6 4  6 6  6 8  6 10  6 12
   7 1  7 3  7 5  7 7  7 9  7 11
   8 0  8 2  8 4  8 6  8 8  8 10  8 12
   9 1  9 3  9 5  9 7  9 9  9 11
   10 0  10 2  10 4  10 6  10 8  10 10 /
\plot 0 0  10 10  8 12  0 4  4 0  10 6  5 11  
   0 6  6 12  10 8  2 0  0 2  9 11 /
\plot  8 0   10 2  1 11  0 10  10 0  /
\plot 0 6  6 0  10 4  3 11  0 8  8 0 /
\setdashes <1mm>
\plot 0 0  0 10 /
\plot 10 0  10 10 /
\setdots <1mm> 
\plot 0 0  10 0 /
%\multiput{$\bullet$} at 0 4  1 9  3 1  3 7  3 11  5 9  6 4  8 2  8 6  8 12  10 4 /
\put{$\ssize P(1)$} at 8.3 13
\put{$\ssize P(4)$} at 2.8 12
\put{$\ssize 1$} at 10 -1
\put{$\ssize 2$} at 8 -1
\put{$\ssize 3$} at 6 -1
\put{$\ssize 4$} at 4 -1
\put{$\ssize 5$} at 2 -1
\put{$\ssize 1$} at 0 -1
\setshadegrid span <.4mm>
\vshade -.2  -.9  2.5  <,z,,> 9.15  8.5  11.9 <z,,,> 11 10  10 /
\vshade 7 0 0  10.2 0 3 /
\vshade -.2 9 10.7 <,z,,> 1.2 10.5  11.8 <z,,,> 2 11 11 /
\vshade -.2 5 6.7 <,z,,> 6.2 11.5  12.8 <z,,,> 7 12 12 /
\vshade 3.5 -.2 .2 <,z,,> 4.3 -.3 0.8  <z,,,> 10.3 5.6 6.9 /

\setshadegrid span <.6mm>
\vshade -.2 5.4 8.8 <,z,,> 5.9 -.3 2.7 <z,,,> 8.5 -.2 0 /
\vshade -.2 1.3 2.7 <,z,,> 1.9 -.3 0.7 <z,,,> 2.6 -.2 0.2 /
\vshade 1 10.2 11.7  10.2 1.3 2.3 / 
\vshade 4.3 11 11 <,z,,> 6.2 9.2  12.7  <z,,,> 10.3 5.2  8.7 /
\endpicture} at 16 0 
%============================================
\put{\beginpicture
\setcoordinatesystem units <.6cm,.6cm>
\multiput{$\bullet$} at 
   0 0  2 0  4 0 
   1 1  3 1
   0 2  2 2  4 2 
   1 3  3 3
   0 4  2 4  4 4 
   /
\plot 0 0  4 4  /
\plot 0 4  4 0 /
\plot 2 0  4 2  2 4  0 2  2 0 /
\setdashes <1mm>
\plot 0 0  0 4 /
\plot 4 0  4 4 /
\setdots <1mm> 
\plot 0 0  4 0 /
\put{$\ssize P(1)$} at  2 4.5
\multiput{$\ssize P(4)$} at 0 4.5  4 4.5 /
\put{$\ssize E(1)$} at 2 -.5
\multiput{$\ssize E(4)$} at 0 -.5  4 -.5 /
\endpicture} at 32 -.4
\put{$\mo \Lambda$} at 8 8
\put{$\Cal C(\Lambda)$} at 32 8

\endpicture}
$$

The resolution quiver $R(\Lambda)$ has the following form:
$$
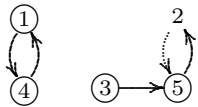
\hbox{\beginpicture
\setcoordinatesystem units <1cm,1cm>
%================================14solid
\put{\beginpicture
\setcoordinatesystem units <.8cm,.8cm>
\put{$\ssize 4$} at 1 0 
\put{$\ssize 1$} at 1 1
\multiput{$\bigcirc$} at 1 0  1 1 /
\arr{0.89 0.22}{0.9 0.2}
\arr{1.11 0.78}{1.1 0.8}
\setquadratic
\plot 1.1 0.2 1.2 0.5 1.1 0.8 /
\plot 0.89 0.2 0.8 0.5 0.9 0.8 /

\endpicture} at 0 0

%================================25halfdotted
\put{\beginpicture
\setcoordinatesystem units <.8cm,.8cm>
\put{$\ssize 5$} at 1 0 
\put{$\ssize 2$} at 1 1
\put{$\ssize 3$} at 0 0
\arr{0.2 0}{0.8 0}
\multiput{$\bigcirc$} at 1 0  0 0  /
\arr{0.89 0.22}{0.9 0.2}
\arr{1.11 0.78}{1.1 0.8}
\setquadratic
\plot 1.1 0.2 1.2 0.5 1.1 0.8 /
\setdots <.4mm>
\plot 0.89 0.2 0.8 0.5 0.9 0.8 /
\endpicture} at 1.3 0
\endpicture}
$$
The black vertices $1,3,4,5$ have been encircled, the arrow $2 \to 5$ has been dotted in
order to stress that it starts at a red vertex. We see that there are two cycles,
one containing only black vertices, the other one containing one black and one red
vertex. The modules in $\Cal X(\Lambda)$ are precisely the simple modules
$S$ which belong to a cycle of black vertices, 
thus we see that $\Cal X(\Lambda) = \{1,4\}$
(as we have mentioned already).

    \bigskip
{\bf Acknowledgment.} This investigation is part of an ongoing 
project with Zhang Pu devoted to determining explicitly the Gorenstein projective modules
for suitable classes of rings and is supported by grant No.~11271251 of NSFC.
The author is grateful to Zhang Pu and his
collaborators at SJTU, but also to Xiao-Wu Chen and the referee for advice and 
very helpful comments, and he has to thank
Julian K\"ulshammer for a careful reading of two versions of the paper.

       \bigskip\bigskip 
%========================================================
{\bf 1.~Notation.}
     \medskip
We denote by $Q = Q(\Lambda)$ the quiver of $\Lambda$; its
vertices are the (isomorphism classes of the) simple $\Lambda$-modules and there
is an arrow $S \to S'$ provided there is a length 2 module with top $S$ and socle
$S'$. Since we assume that no simple module is projective, the quiver $Q$
is just a cycle. 
Let $\tau$ be the Auslander-Reiten translation. Since $\Lambda$ is a Nakayama algebra without simple projective modules, an arrow $S \to S'$ in $Q$
corresponds to the assertion $\tau(S) = S'.$
Often, we will use the vertices $x$ of the quiver $Q$ in order to
index the corresponding modules: thus, we write $S(x)$ for the simple module itself, 
$P(x)$ for the projective cover of $S(x)$.
       \medskip
Given a module $M$, recall that $P(M)$ denotes the projective cover, and we denote by $\Omega(M)$
the first syzygy module. Inductively, let $P_0(M) = P(M)$, $\Omega_0(M) = M$,
and for $n\ge 1$, let
$P_n(M) = P(\Omega_n(M))$ and $\Omega_{n}(M) = \Omega(\Omega_{n-1}(M)).$ Thus, 
there are exact sequences
$$
 0 \to \Omega_{n+1}(M) \to P_n(M) \to \Omega_n(M) \to 0
$$
for all $n\ge 0$. A minimal projective resolution of $M$ has the form
$$
 \cdots \to P_{n}(M) \to P_{n-1}(M) \to \cdots \to P_1(M) \to P_0(M) \to M \to 0,
$$
with $\Omega_n(M)$ being the image of $P_{n}(M) \to P_{n-1}(M).$
     \medskip
Recall that a complex $P_\bullet = (P_i,\delta_i)$ with maps $\delta:P_i \to P_{i-1}$
is called a {\it complete 
projective resolution} provided it is an exact complex of projective modules $P_i$
such that also the complex $\Hom_\Lambda(P_\bullet,\Lambda)$ is exact. The latter
condition is equivalent to the requirement that the inclusion map $\Im(\delta_i)
\subseteq P_{i-1}$ is a left $\Lambda$-approximation, for each $i$. A module $M$
is said to be {\it Gorenstein projective} provided there is a complete
projective resolution $(P_i,\delta_i)$ such that $M = \Im(\delta_0)$. 

An artin algebra $\Lambda$ is said to be a {\it Gorenstein algebra} provided
the injective dimension of ${}_\Lambda\Lambda$ as well as of $\Lambda_\Lambda$ 
is finite. If this is the case, these dimensions are equal and called
the {\it Gorenstein dimension} (or also the {\it virtual dimension}) of $\Lambda$;
we will denote it by $v(\Lambda).$ Note that $\Lambda$ is a Gorenstein algebra
with Gorenstein dimension at most $v$ if and only if $\Omega_v(M)$ is a 
Gorenstein projective module, for every module $M$.

	   \bigskip\bigskip 
%====================================================
{\bf 2.~Projective resolutions of indecomposable $\Lambda$-modules.}
     \medskip
The aim of this section is to point out the relevance of the minimal projective
modules. In particular, we will provide the proof of Proposition 4.
Recall that an indecomposable projective module $P$ is said to be {\it minimal projective} provided
its radical is non-projective (thus provided $P = P(S)$ for some simple module $S$
of projective dimension at least 2). Note that a proper non-zero submodule of
a minimal projective module is not projective (this explains the name). Recall that
a module $M$ is said to be {\it torsionless} provided it is a submodule of a projective
module. 
	\medskip
{\bf Lemma 1.} {\it Let $M$ be an indecomposable non-projective module. If $M$ is torsionless,
then $M$ can be embedded into a minimal projective module and  
any such embedding is a minimal left $\Lambda$-approximation.}
    \medskip
Proof. If there exists an embedding $M \to P$, with $P$ projective, then
there is such an embedding  $M \to P_0$ with $P_0$ indecomposable (since the socle
of $M$ is simple). Since we assume
that $M$ is not projective, there is such an embedding $M\to P_0$, where
$P_0$ is in addition minimal projective. This shows that any torsionless module
can be embedded into a minimal projective module.

Let us fix an embedding $\iota\:M \to P_0$ with $P_0$ minimal projective. In order
to see that $\iota$ is a left $\Lambda$-approximation, we have to show that any non-zero
map $f\:M \to P_1$ with $P_1$ indecomposable projective factors through $\iota$. 
Here, we can assume that also
$P_1$ is minimal projective (namely, since $M$ is not projective, it will map
into $\rad P_1$, thus, if $\rad P_1$ is projective, we may replace $P_1$ by
$\rad P_1$, and so on).

We want to show that $f$ factors through $\iota$. One possibility
is to look at the Auslander-Reiten quiver of $\Lambda$ and consider various
right almost split maps. Here is a proof which uses the fact 
that any automorphism of a submodule of an indecomposable module $X$ 
can be extended to an automorphism of $X$.

First, consider the case where $f$ is a monomorphism.
Thus $P_0$ is a submodule of $P_1$ and the image
of $f$ coincides with the image of $\iota$. It follows that $f$ factors through $\iota$. 

Thus, we can assume that $f$ is not a monomorphism, thus $\Ker(f) \neq 0.$ 
Let $S$ be the socle of the image of $f$ and
take the composition series
$$
 \Ker(f) = K_0 \supset K_1 \supset \cdots \supset K_t = 0
$$
of $\Ker(f)$. Note that $t\ge 1$ and $K_{i-1}/K_i = \tau^i S$ for $1\le i \le t.$

Now $f$ maps into $N = \rad P_1$, and by assumption $N$
is not projective. Let $q\:P = P(N) \to N$ be the projective cover of $N$.
Note that $P(N)$ is both projective and injective. 
Let $U$ be the kernel of $q$ and 
$$
 U = U_0 \supset U_1 \supset \cdots \supset U_u = 0
$$
the composition series of $U$. Then also $U_{i-1}/U_i = \tau^i S$
for $1 \le i \le u.$ For $1 \le i \le u-1$, the module $P/U_i$
is injective and not projective, thus the modules 
$\tau^iS = U_{i-1}/U_i$ for $1\le i \le u-1$
are not torsionless. 

Now assume that $t < u$. Then $1 \le t \le u-1$ and therefore 
$K_{t-1}/K_t = \tau^t S$ is not torsionless. But this is the socle of the
module $M$ and by assumption there is the embedding $\iota\:M \to P_0.$ 
This contradiction shows that we must have $u\le t.$ 

Let $M'$ be the image of $f$ and write $f = f_3f_2f_1$, where $f_1\:M \to M'$, 
whereas $f_2\:M' \to N$
and $f_3\:N \to P_1$ are the inclusion maps. Let $M'' = q^{-1}(M')$ with
inclusion map $f'_2\:M'' \to P(N)$ and let $q'\:M'' \to M'$ be the restriction of $q$,
thus $qf'_2 = f_2q'.$ The length of $M''$ is 
$$
 |M''| = |M'| + |\Ker(q)| = |M'|+u \le |M'|+t  = |M|,
$$
and therefore we can lift the map $f_1\:M \to M'$ to $M''$, thus there is
$f_1'\:M \to M''$ with $q'f_1' = f_1$ and therefore $f_2f_1  = f_2q'f_1' = qf_2'f_1'.$

Since $P(N)$ is injective, there is a map $g\:P_0 \to P(N)$ such
that $g\iota = f_2'f_1'.$ Thus $f = f_3f_2f_1 = f_3qf_2'f_1' = f_3qg\iota$ shows that $f$
factors through $\iota.$
	  \bigskip
Since minimal left approximations are uniquely determined up to isomorphisms, it
follows that {\it an embedding $M \to P$ with $M$ indecomposable and 
$P$ projective is a minimal left 
$\Lambda$-approximation if and only if $P$ is a minimal projective
module.}
	\bigskip
{\bf Proof of Proposition 4.} (i) $\implies$ (ii). This follows from Lemma 1.

(ii) $\implies$ (iii). We assume that we deal with the projective resolution
$$
 \cdots \to P_n \to P_{n-1} \to \cdots \to P_1 \to P_0 \to M \to 0.
$$
Now all the images $\Omega_n(M)$ are indecomposable, thus there are
natural numbers $0 \le n' < n$ such that $\Omega_n(M) = \Omega_{n'}(M)$.
Choose $n$ minimal and assume that $n' \ge 1.$ Now $\Omega_n(M)$ is
a submodule of $P_{n-1}$, and $\Omega_{n'}(M)$ is a submodule of $P_{n'-1}$.
Since $\Omega_n(M) = \Omega_{n'}(M)$, the minimality of $P_{n-1}$ and $P_{n'-1}$
implies that $P_{n-1} = P_{n'-1}$ and therefore $\Omega_{n-1}(M) = \Omega_{n'-1}(M)$.
But this contradicts the minimality of $n$. Thus $n' = 0$ and 
$\Omega_n(M) = \Omega_0(M) = M.$ 

(iii) $\implies$ (i). We assume that there is given an exact sequence
$$
 0 \to M \to P_{n-1} \to \cdots \to P_0 \to M \to 0
$$
such that all the modules $P_i$ are minimal projective modules.
Concatenation of countably many such sequences yields a complete resolution
with $M$ as one of the images. This shows that $M$ is Gorenstein projective.
     \bigskip\bigskip 
%================================================================
{\bf 3.~The resolution quiver of $\Lambda$.}
     \bigskip 
Let $\gamma S = \tau\soc P(S)$.
The importance of the map $\gamma$ stems from two observations, see Lemmas 2 and 3.
As Xiao-Wu Chen has pointed out, the map $\gamma$ has been
considered already by W. H. Gustafson [G]. 
	   \medskip
{\bf Lemma 2.} {\it Let $M$ be an indecomposable module.
Then either the projective dimension of $M$ is at most 1 and $\Omega_2(M) = 0$, 
or else $\top \Omega_2(M) = \gamma \top M.$}
   \medskip
Proof: Write $M = P(M)/U$ for some submodule $U$ of $P(M)$. 
We can assume that $U$ is a proper submodule of $P(M)$. 
There are exact sequences of the following form
$$
\gather
 0 \to \Omega_2(M) \to P_1(M) \to \Omega_1(M) \to 0 \cr
 0 \to \Omega_1(M) \to P(M) \to M \to 0.
\endgather
$$
Clearly, $\soc \Omega_1(M) = \soc P(M)$ and the first exact sequence
shows that either $\Omega_2(M) = 0$ (thus the projective dimension of $M$
is at most 1)
or else $\top \Omega_2(M) = \tau \soc \Omega_1(M)$. In the latter case, 
$\top \Omega_2(M) = \gamma \top M$.
      \medskip
Inductively, we see: 
	     \medskip
{\bf Corollary.} {\it Let $M$ be an indecomposable module and $m\in \Bbb N$. 
Then either $\Omega_{2m}(M) = 0$ or else $\top \Omega_{2m}(M)  = \gamma^m\top M.$}
     \bigskip\bigskip
%================================================	
If $|P(S)| \ge s$, we define $H(S)$ to be the factor module of $P(S)$ of length 
$s$, where $S$ is any simple module. We call this module $H(S)$ a {\it primitive} module.  
In case $p > s$, these modules $H(S)$ do exist for all simple modules $S$ and
are non-projective. 
    \medskip
{\bf Lemma 3.} {\it Let $S$ be a simple module with $|P(S)| > s$.
Then $|P(\gamma S)| \ge s$ and}
$$
 \Omega_2 H(S) = H(\gamma S).
$$
	\medskip
Proof. Since $H(S)$ has length $s$, and $|P(S)| > s$, 
the minimal projective presentation of $H(S)$ is
$$
 P(S) \to P(S) \to H(S) \to 0,
$$
and therefore $\Omega_2 H(S)$ has length $s$. On the other hand, according
to Lemma 2, we know that $\Omega_2 H(S)$ is a factor module of $P(\gamma S)$.
This yields both assertions.
     \medskip

{\bf Corollary.} {\it Let $p>s$. Let $S$ be a simple module and $m$ a natural number. Then} 
$$
 \Omega_{2m} H(S) = H(\gamma^m S).
$$
	\medskip
Note that this corollary implies that the projective dimension of $H(S)$
is infinite. Thus, any Nakayama algebra with $p>s$ has infinite global dimension.
This implies the following result of Gustafsen [G]: if the Loewy length of
$\Lambda$ is greater than or equal to $2s$, then $\Lambda$ has infinite global
dimension (since in this case $p > s$). 
	  \bigskip 
The map $\gamma$ is used in order to 
obtain the {\it resolution quiver} $R = R(\Lambda)$ (as introduced in the
introduction): its vertices are the (isomorphism classes of the)
simple $\Lambda$-modules, and for any simple module $S$, 
there is an arrow $S \to \gamma S$.

Note that {\it any connected component of the resolution quiver has a unique cycle.}
This follows from the fact that at any vertex $x$ precisely one arrow
starts; thus given any connected component, the number of arrows in the component
is equal to the number of vertices in the component. 
   \medskip
We say that a vertex $x$ or the corresponding simple or projective modules $S(x)$
and $P(x)$ are {\it black} provided $P(x)$ is minimal projective, otherwise $x$ (and
$S(x)$ and $P(x)$) will be said to be {\it red}. Note that $x$ is red if and only if
the projective dimension of $S(x)$ is equal to $1$, and black if and only if 
the projective dimension of $S(x)$ is greater than or equal to $2$.
A cycle in $R$ 
will be said to be {\it black} provided all the vertices occurring
in the cycle are black. A vertex $x$
is said to be {\it cyclically black} provided $x$ belongs to a black cycle.
   \bigskip
The resolution quiver has the following property: {\it any vertex $y$ is end point
of at most one arrow $x\to y$ with $x$ black.} Namely, if $S,S'$ are black
simple modules with $\gamma S = \gamma S'$, then $\tau\soc P(S) = \tau\soc P(S')$,
thus $\soc P(S) = \soc P(S')$ and therefore $P(S) \subseteq P(S')$
or $P(S') \subseteq P(S)$. However, since both $P(S)$ and $P(S')$ are
minimal projective modules, it follows that $P(S) = P(S')$ and thus $S = S'.$ 
	\medskip
As a consequence we obtain:
   \medskip
{\bf Lemma 4 (Red Entrance Lemma).}
{\it If $x_0 \to x_1 \to \cdots \to x_a \to y$ is a path 
such that $y$ is cyclically black, whereas $x_{a}$ is not cyclically black,
then $x_{a}$ is red.} 
     \medskip
Proof. Since $y$ is cyclically black, there is an arrow
$x' \to y$ such that also $x'$ is cyclically black. Since $x_a$ is not
cyclically black, we have $x_{a} \neq x'$, and it 
follows that $x_{a}$ cannot be black, thus it is red. 
	\medskip
There is the following consequence: {\it Let $x$ be a vertex which is 
not cyclically black. Then there exists $m\ge 0$ such that $\gamma^mx$ is red.}
    \bigskip
{\bf Lemma 5.} 
{\it Let $M$ be indecomposable and not projective.  Then the following conditions
are equivalent:}
\item{\rm(i)} {\it $M$ is Gorenstein projective.} 
\item{\rm(ii)} {\it Both $\top M$ and $\top \Omega(M)$
are cyclically black.}
\item{\rm(iii)} {\it Both $\top M$ and $\tau\soc M$
are cyclically black.}
    \medskip
Proof: 
The assertions (ii) and (iii) are equivalent, since $\tau\soc M = \top\Omega(M).$
    \medskip
(ii) $\implies$ (i). If $\top M$ is cyclically black, then $M$ has infinite projective dimension 
and all the modules $P_{2m}(M)$ are black, for $m\ge 0$. 
If $\top \Omega(M)$ is cyclically black, then all the modules 
$P_{2m+1}(M)$ are black, for $m\ge 0$. Thus, if both conditions are satisfied, then
all the projective modules occurring in the minimal resolution of $M$ are
minimal, thus $M$ is Gorenstein projective, according to Lemma 2.

(i) $\implies$ (ii). Assume that $\top M$ is not cyclically black. Then there is 
some $m\ge 0$ such that $\gamma^m\top M$ is red, but $P_{2m}(M) = P(\gamma^m \top M)$.
Similarly, if  $\top \Omega M$ is not cyclically black, then there is 
some $m\ge 0$ such that $\gamma^m\top\Omega  M$ is red, and 
$P_{2m+1}(M) = P(\gamma^m \top\Omega M)$. Thus, $M$ cannot be Gorenstein projective,
according to Lemma 2.
	  \medskip
Several consequences are of interest. First, we look at part (b) of Proposition 5. We see that
one direction works without the assumption $p > s$.
    \medskip 
{\bf Corollary 1.} {\it 
Let $\Lambda$ be a Nakayama algebra and assume that no cycle of $R(\Lambda)$ is  black. Then $\Lambda$ is $\CM$-free.}
    \bigskip
{\bf Corollary 2.} {\it If $S$ is a cyclically black simple module
and $|P(S)| > s$, then $H(S)$ is Gorenstein projective and not projective.}
    \medskip
Proof: Let $S$ by cyclically black and $|P(S)| > s$. Then $H(S)$
exists and is non-projective. Both $\top H(S)$ and $\top \Omega H(S)$ are
equal to $S$, thus cyclically black. According to Lemma 5, we see that 
$H(S)$ is Gorenstein projective. 
       \bigskip
{\bf Corollary 3.} {\it The simple modules in $\Cal X(\Lambda)$ are
cyclically black. Conversely, if $S$ is a cyclically black simple module
and $|P(S)| > s$, then $S$ belongs to
$\Cal X(\Lambda).$}
      \medskip
Proof. Let $S$ be in $\Cal X(\Lambda).$ Then $P(S)$ has a factor module $M$
which is Gorenstein projective, but not projective. Lemma 5 asserts that
$S = \top M$ is cyclically black.

On the other hand, let $S$ by cyclically black and $|P(S)| > s$. Then $H(S)$
exists and is non-projective. According to Corollary 2, $H(S)$ is
Gorenstein projective. Since $P(S)$ has a non-projective
factor module which is Gorenstein projective, we see that $S$ belongs to
$\Cal X(\Lambda).$
      \medskip
{\bf Proof of Proposition 5 (b).} One direction is covered by Corollary 1.
For the converse, let us assume that $p > s$ and that there is a black cycle,
say containing the vertex $x$. 
According to Corollary 2, the module $H(x)$ is Gorenstein projective and not
projective, thus $\Lambda$ is not $\CM$-free.
	    \medskip 
{\bf Proof of Proposition 3.}
By definition, $\Cal X = \Cal X(\Lambda)$ is the class of simple modules $S$
such that $P(S)$ has a factor module which is non-projective and
Gorenstein projective, thus $S$ belongs to $\Cal X$ if and only if
there is a non-projective Gorenstein projective module $N$ with $S = \top N.$

To show the equivalence of (i), (ii), (iii), it is sufficient to consider the
case of $M$ being indecomposable. 

First, assume that $M$ is projective, say $M = P(S)$ for some simple module $S$.
If $P(S) \in \Cal C$, then $P(S)$ is the projective cover of a non-projective
Gorenstein projective module $M'$, thus 
$S\in \Cal X$, this shows that (i) implies (iii). In order to show that (i)
implies (ii), we have to show in addition that $\tau\soc M$ belongs to $\Cal X.$ 
Since $M'$ is non-projective and Gorenstein projective, also $\Omega^2 M'$ is
non-projective and Gorenstein projective, thus $\tau \soc M = \tau\soc \Omega M' =
\top \Omega^2 M'$ belongs to $\Cal X$.
Conversely, if (ii) or (iii) is satisfied, then
$\top M$ belongs to $\Cal X$, thus $M = P(\top M)$ belongs to $\Cal C$.

Next, assume that $M$ is non-projective. The assertions (ii) and (iii)
are equivalent since for $M$ indecomposable and non-projective, $\tau\soc M =
\top\Omega M.$ If $M$ is Gorenstein projective, then also $\Omega M$ is
Gorenstein projective (and non-projective), thus $\top M$ and $ \top \Omega M$
both belong to $\Cal X.$ This shows that (i) implies (iii). Conversely,
assume that $\top M$ and $\top \Omega M$ are in $\Cal X$. Then, according to
the Corollary 3, $\top M$ and $\top \Omega M$ are cyclically black
and therefore Lemma 5 asserts that $M$ is Gorenstein projective. This shows that
(iii) implies (i) and completes the proof of Proposition 3. 
      \bigskip 
{\bf Remark 1.} It seems to be of interest to identify the sources of the
resolution quiver $R$: {\it The simple module $S$ is a source of $R$ 
if and only if $\tau^-S$ cannot be embedded into $\Lambda$} (thus if and only
if $I(\tau^-S)$ is not projective). Namely, 
if $S$ is in the image of $\gamma$, then there
is a simple module $S'$ with $S = \gamma S' = \tau\soc P(S')$, thus $\tau^-S =
\soc P(S') \subseteq \soc \Lambda$. And conversely, if $\tau^-S$ 
can be embedded into $\Lambda$, then
it can be embedded into some indecomposable projective module $P(S')$, and then
$\gamma S' = \tau\soc P(S') = \tau\tau^-S = S,$  thus $S$ is in the image of $\gamma.$

As a consequence, {\it the number of sources of $R$ is equal to $s-t$,} where $t$ is 
the
number of indecomposable modules which are both projective and injective (and this
is also the number of minimal projective modules).

The same argument shows that in general {\it the number of arrows ending in $S$
is equal to the number of projective modules with socle $\tau^-S.$}
   \bigskip
{\bf Remark 2.} The referee has pointed out that a simple module $S$ is cyclically black if and only if it is perfect in the sense of [CY]. 
     \bigskip\bigskip 

%==============================
{\bf 4.~The elementary Gorenstein projective modules.}
     \medskip
We are going to use Lemma 5 (or else Proposition 3) in order to show
some important 
closure properties of $\gp_0$. We need them in order to prove Proposition 1.
After Proposition 1 is established, they are direct consequences. 
      \medskip 
{\bf Lemma 6.} {\it Let $X,Y$ be indecomposable modules which are
Gorenstein projective and not projective. If $f\:X \to Y$ is a homomorphism,
then kernel, image and cokernel of $f$ are Gorenstein projective modules.}
     \medskip 
Proof: We can assume that $f$ is non-zero. 
Then the image $Z$ of $f$ is indecomposable and not projective.
According to Lemma 5, $\top X, \top Y, \tau\soc X, \tau\soc Y$ all 
are cyclically black. 
But $\top Z = \top X$ and $\soc Z = \soc Y$, thus $\top Z$ and
$\tau\soc Z$ are cyclically black. Using again Lemma 5, we see that
$Z$ is Gorenstein projective.
If the kernel $K$ of $f$ is non-zero, then $K$ is indecomposable and non-projective.
Also, $\soc K = \soc X$, and 
$\top K = \tau\soc Z = \tau \soc Y$, thus Lemma 5 implies again that $K$ is
Gorenstein projective.  Finally, if the cokernel $Q$ of $f$ is non-zero, 
it is indecomposable and non-projective, and $\top Q = \top Y$, and
$\tau\soc Q = \top Z = \top X$ are cyclically black, so that also $Q$
is Gorenstein projective. 
   \bigskip 
{\bf Lemma 7.} {\it Let $X$ be an indecomposable module which is Gorenstein
projective and not projective. If $X' \subseteq X''$ are submodules of $X$
such that $X''/X'$ is Gorenstein projective and non-zero, then also
$X'$ and $X/X''$ are Gorenstein projective.}
     \medskip
Proof: Again, we use Lemma 5. 
Of course, also $X''/X'$ is non-projective, thus 
$\top X, \tau\soc X, \top (X''/X'), \tau\soc(X''/X')$ are 
cyclically black. If $X' \neq 0$, then this is a non-projective indecomposable
module with $\top X' = \tau\soc(X''/X')$ and $\tau\soc X' = \tau\soc X$
cyclically black. If $X/X'' \neq 0$, then this is a non-projective indecomposable
module with $\top(X/X'') = \top X$ and $\tau\soc(X/X'') = \top (X''/X')$
cyclically black.
	   \bigskip

Let $S$ be a simple module. If $P(S)$ has a non-projective factor module
which is Gorenstein projective, then we denote the smallest module of this kind
by $E(S).$ If $M$ is any module, then $M$ is of the form $E(S)$ 
provided $M$ is non-projective, but Gorenstein projective, and the only 
proper factor module of $M$ which is Gorenstein projective, is the zero module
(and $S = \top M$). Thus, the modules $E(S)$ are the elementary Gorenstein 
projective modules as defined in the introduction. 
	   \medskip
Let $g$ be the number of isomorphism classes of elementary Gorenstein projective modules.
Note that $g$ equals the number of simple modules in $\Cal X(\Lambda).$ 
    \medskip
{\bf Proof of Proposition 2.} 
Let $E, E'$ be elementary modules and $f\:E \to E'$ a non-zero homomorphism.
The image of $f$ is a non-projective
Gorenstein projective factor module of $E$, thus $f$ has to be injective,
according to Lemma 6.
If $f$ is not surjective, then the cokernel of $f$ is a non-projective
Gorenstein projective proper factor module of $E'$, impossible (again we use
Lemma 6).  This shows that
$f$ is an isomorphism. Thus, if $E_1,\dots, E_g$ are representatives
of the isomorphism classes of $\Cal E$, then this is a set of pairwise
orthogonal bricks. 

If $U$ is a proper submodule of $E$ which is Gorenstein projective,
then $E/U$ is a non-zero factor module of $E$ which is Gorenstein projective,
according to Lemma 7. But this implies that $U = 0$. This shows that
$\Cal E \subseteq \Cal E'.$ 

Conversely, assume that $X$ belongs to $\Cal E'$, thus $X$ is an indecomposable
module which is Gorenstein projective and non-projective, such that
no proper non-zero submodule is Gorenstein projective. Assume that there is
a proper factor module $X/U$ of $X$ which is Gorenstein projective. Then
Lemma 6 asserts that $U$ is Gorenstein projective and therefore $U = X$.
This shows that $\Cal E' \subseteq \Cal E.$ Thus we have established (c).

Next, let us show that $\Cal C = \Cal F(\Cal E)$. Since we deal with a set
of orthogonal bricks, $\Cal F(\Cal E)$ is an abelian subcategory. 
Since all the modules in $\Cal E$ are Gorenstein projective, the modules
in $\Cal F(\Cal E)$ are Gorenstein projective. Assume that $P$ is indecomposable
projective and in $\Cal F(\Cal E)$. Then $P$ has a non-zero submodule $U$ such that
$P/U$ belongs to $\Cal E$. But this shows that $P$ is the projective cover of 
a Gorenstein projective module which is not projective.
This shows that $\Cal F(\Cal E) \subseteq \Cal C$. 

Conversely, consider
an indecomposable module $X$ in $\Cal C$. For the non-projective
modules we use induction on the length. Thus, 
let $X$ be non-projective, let $P(X) = P(S)$ be its projective cover, with
$S$ simple. Since $X$ is Gorenstein projective, the module $E(S)$ exists and
is a factor module of $X$, say $X/U = E(S)$. Now $U$ itself is Gorenstein
projective, according to Lemma 6, thus by induction it belongs to $\Cal F(\Cal E).$
Then also $X$ belongs to $\Cal F(\Cal E).$ 

If $X = P(S)$ is projective, then $X$
is the projective cover of a Gorenstein projective module which is not projective,
thus the factor module $E(S)$ exists, say $E(S) = X/X'$. But $X'$
is Gorenstein projective and not projective, thus we know already that
$X'$ belongs to $\Cal F(\Cal E).$ This completes the proof of (a).
     \medskip
It remains to establish (b). We assume now that $\Lambda$ is connected
and $\Cal C$ is non-zero, and we denote by $E_1,\dots, E_g$ representatives 
of the isomorphism classes in $\Cal E.$ On the one hand,
these modules are pairwise orthogonal bricks. On the other hand, according to
Lemma 7, no module $E_i$ is a proper subquotient of some $E_j$, thus the 
modules $E_i$ must have pairwise different support. Thus, any simple module
occurs with multiplicity at most 1 in $\bigoplus_{i=1}^g E_i.$

It remains to be seen that any simple module occurs in the support of some
$E_i$. Since $\Cal C$ is not zero, we have $g \ge 1$, thus there is a simple
module $S$ which belongs to the support of $\bigoplus_{i=1}^g E_i.$ We show:
if $S$ is a composition factor of an elementary module 
$E$, then $\tau S$ is a composition factor of  an elementary module $E'$.
Thus assume that $S$ is a composition
factor of $E$. If $S$ is not the socle of $E$, then $\tau S$ is also
composition factor of $E$. Thus, we have to consider the case that $S = \soc E$. 
The following lemma shows that there is an elementary module $E'$ such that
$\tau S$ is a composition factor of $E'$. 
Since we assume that $\Lambda$ is connected,
the set of simple $\Lambda$-modules is a single $\tau$-orbit, therefore
any simple module occurs as a composition factor of some elementary module.
This completes the proof of (b).
     \bigskip 
     
{\bf Lemma 8.} {\it If $E$ is an elementary module, then there is an elementary
module $E'$ with $\top E' = \tau\soc E$ and $\Ext^1(E,E') \neq 0.$}
       \medskip 
Proof. Let $P(E)$
be the projective cover of $E$. This is a module in $\Cal C = \Cal F(\Cal E)$,
thus $P(E)$ has a filtration using elementary Gorenstein projective modules,
say
$$
 P(E) = M_0 \supset M_1 \supset \cdots \supset M_t = 0,
$$
such that $M_{i-1}/M_i$ is elementary, for all $1\le i \le t.$ 
Note that $M_0/M_1 = E$. Since $E$
is not projective, we have $t\ge 2$, thus let $E' = M_1/M_2.$ Then
$$
 \top E' = \top (M_1/M_2) = \tau \soc (M_0/M_1) = \tau\soc E.
$$ 
Since the module $M_0/M_2$ is indecomposable,
the exact sequence
$$
 0 \to M_1/M_2 \to M_0/M_2 \to M_0/M_1 \to 0
$$
does not split, therefore $\Ext^1(E,E') \neq 0.$
     \bigskip

{\bf Proof of Proposition 1.}
Since $\Cal C = \Cal F(\Cal E),$ and $\Cal E$ is a set of pairwise orthogonal
bricks, $\Cal C$ is a full exact subcategory of $\mo \Lambda$ which is
closed under extensions. By definition, the category $\Cal C$ is closed under
projective covers. In order to see that $\Cal C$ is closed under
minimal left $\Lambda$-approximations, consider an indecomposable object $M$
in $\Cal C$. We may assume that $M$ is not projective. According to Proposition 4,
there is an exact sequence
$$
 0 \to M \to P_{n-1} \to \cdots \to P_0 \to M \to 0
$$
such that all the modules $P_i$ are minimal projective modules. With $M$
also $\Omega_{n-1}(M)$ is Gorenstein projective and not projective, and
$P_{n-1}$ is its projective cover. By the construction of $\Cal C$, the module
$P_{n-1}$ belongs to $\Cal C$. Lemma 1 asserts that the embedding
$M \to P_{n-1}$ is a minimal left $\Lambda$-approximation.

As an abelian length category with only finitely many indecomposable objects,
the category $\Cal C$ is equivalent to the module category of an  
algebra $\Lambda'$, and $\Lambda'$ has precisely $g$ simple
modules. The indecomposable objects in $\Cal C$ are indecomposable 
$\Lambda$-modules: since they have a unique composition series as $\Lambda$-modules,
they also have a unique composition series inside the category $\Cal C$. This
shows that $\Lambda'$ is again a Nakayama algebra. The indecomposable 
projective objects in $\Cal C$ are the projective $\Lambda$-modules which belong to
$\Cal C$: these are some of the minimal projective $\Lambda$-modules.
But a minimal projective $\Lambda$-module cannot be properly embedded into
any other minimal projective $\Lambda$-module. This shows that the indecomposable
projective objects in $\Cal C$ are injective in $\Cal C$, thus $\Lambda'$
is self-injective.

If $\Lambda$ is connected and $\Cal C \neq 0$, then we want to see that $\Lambda'$
is connected. But this follows immediately from Proposition 2 and Lemma 8.

It remains to show that the embedding $\mo\Lambda' \to \mo\Lambda$ induces an embedding
$\underline{\mo}\Lambda' \to \underline{\mo}\Lambda$ with image just
$\underline{\gp}$.  Let $X,Y$ be indecomposable in $\gp_0$, and $f\:X \to Y$ a
morphism which factors through a projective $\Lambda$-module, say $P$.
Then the morphism $X \to P$ factors through $P'$, where $X\to P'$ is the
minimal left $\Lambda$-approximation of $X$. But $P'$ belongs to $\Cal C$
and is projective in $\Cal C$, thus $f$ is zero in $\underline{\gp}.$
    \bigskip

{\bf More information about $\Lambda'$.} We know that
the number of simple $\Lambda'$-modules is $g$. Let us insert a formula
for the length $p'$ of the indecomposable projective $\Lambda'$-modules.
We assume that $\Lambda$ is a connected Nakayama algebra and that $\Cal C$ is non-zero.
Let $E_1,\dots, E_g$ be a complete set of elementary Gorenstein projective
modules and let $p_i = |P(E_i)|.$ Then 
{\it the length $p'$ of the indecomposable projective $\Lambda'$-modules is given by
the formula}
$$
 p' = \frac1s\sum_{i=1}^g p_i
$$
	\medskip
Proof: Let $E_i = E_j$ provided $i\equiv j\mod g$. We can assume that
$\Ext^1(E_i,E_{i+1})\neq 0$ for all $i$, thus $P(E_i)$ has in $\Cal C$ 
a composition series with factors $E_i, E_{i+1},\dots, E_{i+p'-1}$ going down.
Then 
$$
\align
 \sum_{i=1}^g p_i &= \sum_{i=1}^g |P(E_i)| 
  =  \sum_{i=1}^{g}\sum_{j=0}^{p'-1} |E_{i+j}| = 
 \sum_{j=0}^{p'-1}\sum_{i=1}^{g} |E_{i+j}| \cr
 &= 
 \sum_{i=0}^{p'-1} s  = p's.
\endalign
$$
Here, we have used that for any $j$ the sequence $E_{j+1},\dots,E_{j+g}$
is obtained from the sequence $E_1,\dots,E_g$ just by permutation and that 
$\sum_{i=1}^g |E_i| = s$, see Proposition 2.
	      \bigskip\bigskip

%======================================================
{\bf 5.~Gorenstein algebras.}
     \medskip
We want to present a proof of part (a) of Proposition 5. We assume in this section
that $p > s.$
     \medskip 
{\bf Proposition 6.} {\it  Let $p > s$. The Nakayama algebra  $\Lambda$ is a
Gorenstein algebra if and only if 
all cycles in the resolution quiver of $\Lambda$ are black. In this case,
the Gorenstein dimension of $\Lambda$ 
is equal to $2d$, where $d$ is the maximal distance
between vertices and the black cycles.}
	\medskip
We need the following lemmas.
   \medskip
{\bf Lemma 9.} {\it Let $x = x_0 \to \cdots \to  x_d$ be a path in the resolution quiver such 
that $x_{d}$ is cyclically black, whereas $x_{d-1}$ is not. Then
$G$-$\dim H(x) = 2d$.}
	  \medskip
Proof. We show that $\Omega_{2d-1}H(x)$ is not Gorenstein projective,
whereas $\Omega_{2d}H(x)$ is Gorenstein projective. According to Lemma 4
we know that $\Omega_{2d} H(x) = H(\gamma^d x) = H(x_d)$ and this is 
a Gorenstein projective module.
On the other hand, there is the following exact sequence
$$
 0 @>>> H(x_d) @>>> P(x_{d-1})  @>f>> P(x_{d-1})  @>>> H(x_{d-1}) @>>> 0
$$
and $\Omega_{2d-1}H(x)$ is the image of $f$. In particular, we see that
the top of $\Omega_{2d-1}H(x)$ is $S(x_{d-1})$ and this is a red vertex,
according to the red entrance lemma. Since the top of $\Omega_{2d-1}H(x)$
is not even black, $\Omega_{2d-1}H(x)$ cannot be Gorenstein projective.
   \medskip
{\bf Lemma 10.} {\it Assume that for any path $x_0 \to \cdots \to  x_d$ 
in the resolution quiver, the vertex $x_d$ is cyclically black. Then the
$G$-dimension of any $\Lambda$-module is at most $2d$.}
	      \medskip
Proof. We show that $\Omega_{2d}(M)$ is Gorenstein projective, 
for any indecomposable module $M$. According to Lemma 5, we have to show that
the modules $\top \Omega_{2d}(M)$ and $\top \Omega_{2d+1}(M)$ are zero or
cyclically black. But $\top \Omega_{2d}(M)$ is zero or equal to
$\gamma^d \top M$, and $\top \Omega_{2d+1}(M)$ is zero or equal to
$\gamma^d \top \Omega M$. The assumption of the lemma can be rephrased by saying
that $\gamma^d S$ is cyclically black for all simple modules $S$. This completes
the proof.
    \medskip
{\bf Proof of Proposition 6.}
First, assume that there is a  cycle in the resolution quiver which involves a red vertex, say the red vertex $x$. Then a minimal projective resolution
of $H(x)$ involves infinitely many copies of $P(x)$ and $P(x)$ is not minimal
projective, thus $H(x)$ has infinite $G$-dimension. This shows that $\Lambda$
is not a Gorenstein algebra. 

On the other hand, if all the cycles in the resolution quiver are black, then 
all indecomposable modules have finite $G$-dimension, according to Lemma 8.
Thus $\Lambda$ is a Gorenstein algebra. 
     \bigskip
{\bf Proposition 7.} {\it If $\Lambda$ is a Gorenstein algebra with $p > s$, then
$v(\Lambda)\le 2s-2.$}
	       \medskip

Proof. Any path 
$x = x_0 \to \cdots \to  x_d$ in $R$ such that $x_{d-1}$ does not belong
to a cycle involves pairwise different vertices, thus $d \le s-1$. Thus, Lemma 8
asserts that the $G$-dimension of any $\Lambda$-module is at most $2s-2$.
	\medskip
As Xiao-Wu Chen has pointed out, this result (and its proof) corresponds to Gustafson's 
bound for the finitistic dimension of a Nakayama algebra.
      \bigskip\bigskip 
%===========================================================
{\bf 6.~Calculation of resolution quivers.}
     \medskip
Recall that $(p_1,\dots, p_s)$ is said to be a {\it Kupisch series} for $\Lambda$,
provided we have labeled the indecomposable projective modules $P_1,\dots, P_s$
such that $\rad P_i$ is a factor module of $P_{i+1}$ (thus provided there is an
arrow $S_i \to S_{i+1}$) and $p_i = |P_i|.$ The Kupisch series for $\Lambda$
are obtained from any one by cyclic permutation. 

It is an easy exercise to draw the resolution quiver of $\Lambda$ if a Kupisch series
is known. Namely, {\it if $(p_1,\dots,p_s)$ is a Kupisch series for $\Lambda$, 
then }
$$
 \gamma (i) \equiv i+p_i \mod s.
$$
	\medskip
The {\it roof} of the Nakayama algebra $\Lambda$ is the factor algebra
$r(\Lambda) = \Lambda/\soc_{p-2}\Lambda$. (Here, for any left module $M$,
the socle sequence $0 = \soc_0M \subseteq \soc_1M \subseteq \soc_2M \subseteq \cdots\subseteq M$ is defined
inductively by $\soc_iM/\soc_{i-1}M = \soc(M/\soc_{i-1}M)$, for all $i\ge 1;$ note that for 
$M= {}_\Lambda\Lambda,$ all the submodules $\soc_iM$ are two-sided ideals.) It is easy to check that
two Nakayama algebras $\Lambda, \Lambda'$ have the same roof 
provided they have Kupisch series $(p_1,\dots,p_s)$ and $(p'_1,\dots,p'_s)$
such that $p_i-p_i' = c$, for $1\le i \le s,$ where $c = c(\Lambda,\Lambda')$
is a constant. Let us also stress: {\it If $c(\Lambda,\Lambda')$ is a multiple of $s$, then
$\Lambda$ and $\Lambda'$ have the same resolution quiver.}

In order to visualize a roof $r(\Lambda)$, we will draw its Auslander-Reiten quiver. Actually,
instead of drawing the Auslander-Reiten quiver of $r(\Lambda)$, we draw the
Auslander-Reiten quiver $\Gamma$ of a Nakayama algebra with linearly directed quiver 
of type $\Bbb A_{s+1}$
such that the Auslander-Reiten quiver of $r(\Lambda)$ is obtained from $\Gamma$ by identifying the
simple projective with the simple injective module.

For example, 
here is the way we present the roof of the Nakayama algebra $\Lambda$ with Kupisch series $(m+3,m+3,m+2)$, 
note that the roof $r(\Lambda)$ is the Nakayama algebra 
with Kupisch series $(3,3,2)$:
$$
\hbox{\beginpicture
\setcoordinatesystem units <.5cm,.5cm>
\multiput{$\bullet$} at 0 0  1 1   2 0  4 0  6 0  3 1  5 1  4 2  2 2 /
\plot 0 0  2 2  4 0  5 1 /
\plot 4 2  6 0 /
\plot 1 1  2 0  4 2 /
\put{$\bigcirc$} at 1 1
\put{$\bigcirc$} at 4 2
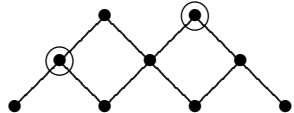
\endpicture} 
$$
This is the Auslander-Reiten quiver
of a Nakayama algebra with linearly directed quiver 
of type $\Bbb A_{4}$, in order to obtain the Auslander-Reiten quiver of $r(\Lambda)$,
the vertex far right
has to be identified with the vertex far left (so that the Auslander-Reiten of $r(\Lambda)$
quiver lives on a cylinder). 

Let us stress the conventions which we use: When we draw an Auslander-Reiten quiver, 
always the direction of the arrows will be deleted: 
all arrows are supposed to be directed from left to right.  The projective vertices of $r(\Lambda)$ (and thus of $\Lambda)$) should be
labeled {\bf going from right to left} as $P(1), P(2), P(3)$, see 
$$
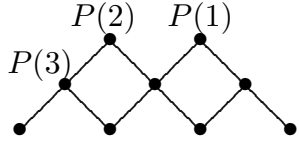
\hbox{\beginpicture
\setcoordinatesystem units <.5cm,.5cm>
\multiput{$\bullet$} at 0 0  1 1   2 0  4 0  6 0  3 1  5 1  4 2  2 2 /
\plot 0 0  2 2  4 0  5 1 /
\plot 4 2  6 0 /
\plot 1 1  2 0  4 2 /
\put{$P(1)$} at 4 2.5
\put{$P(2)$} at 1.85 2.5
\put{$P(3)$} at 0.45 1.45

\endpicture} 
$$
In the previous picture, two of the projective vertices,
namely $P(1)$ and $P(3)$, have been encircled, since
they correspond to the minimal projective modules of $r(\Lambda)$ and of $\Lambda$. 

In order to draw the resolution quiver, we use as vertices the numbers $1,2,\dots, s$
and we draw an arrow $i \to j$ provided $j \equiv i+p_i \mod s.$ 
For the convenience of the reader, the black vertices will be encircled and the 
arrows which start at a red vertex will be dotted (this is quite redundant, but maybe helpful).
Thus, the resolution quiver
for the algebra with Kupisch series $(5,5,4)$ (or, more generally, with Kupisch series 
$(3m+2,3m+2,3m+1)$ for some $m\ge 0$) looks as
follows:
$$
\hbox{\beginpicture
\setcoordinatesystem units <1cm,1cm>
\put{$\ssize 2$} at 0 0
\put{$\ssize 1$} at 1 0
\put{$\ssize 3$} at 2 0
\multiput{$\bigcirc$} at 2 0  1 0 /

\setquadratic
\plot 1.2 0.1  1.5 0.2  1.8 0.1 /
\arr{1.75 0.125}{1.8 0.1}
%\circulararc 300 degrees from 1.1 -.2  center at 1.5 0 
\plot 1.2 -.1  1.5 -.2  1.8 -.1 /
\arr{1.25 -.125}{1.2 -.1}
\setdots <.4mm>
\arr{0.2 0}{0.8 0}

\endpicture} 
$$
	\bigskip
The bound $2s-2$ in Proposition 7 is optimal as the algebras with Kupisch series
$$
 (ms+1,ms+1,\dots,ms+1,ms)
$$
and $m\in \Bbb N_1$ show. Namely, the resolution quiver looks as follows
$$
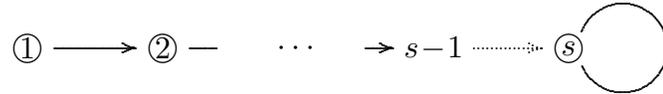
{\beginpicture
\setcoordinatesystem units <1.5cm,1cm>
\multiput{} at 0 0  0 -.2 /
\put{$1$} at 0 0 
\put{$2$} at 1 0 
\put{$\cdots$} at 2 0 
\put{$s\!-\!1$} at 3 0 
\put{$s$} at 4 0 
\arr{0.2 0}{0.8 0}
\plot 1.2 0  1.4 0 /
\arr{2.5 0}{2.7 0}
\circulararc 315 degrees from 4.1 -.2  center at 4.4 0 
\setdots <.5mm>
\arr{3.3 0}{3.8 0}
\multiput{$\bigcirc$} at 0 0  1 0  4 0 /
\endpicture}
$$
(since $\gamma s = s+p_s =  s + ms \equiv s \mod s$, whereas
$\gamma i = i+p_i = i+ms+1 \equiv i+1$ for $1 \le i \le s-1$).
The vertex $s-1$ is red, all others are black. In particular, the
loop at the vertex $s$ is black. Thus, according to Lemma 9, 
the path $1 \to 2 \to \cdots 
\to s\!-\!1 \to s$ shows that the $G$-dimension of $H(s)$ is equal to
$2s-2.$ 
	\bigskip
{\bf Examples.}
We present the different types of the connected Nakayama algebras with $s =3, 4, 5$
(and $p > s$). Note that for $s=3$, 
this classification can be found already in the paper [CY] by Chen and Ye.
     \bigskip
In the tables,  we show on the left
the different roofs (or at least the upper boundary of the roofs).
The remaining columns are indexed 
by the numbers $a$ with $1 \le a \le s$ and show for $s=3,4$ 
the resolution quiver $R(\Lambda)$
of the Nakayama algebras $\Lambda$ with given roof and $p \equiv a \mod s,$
for $s=5$ only the cycles of $R(\Lambda)$ are exhibited.
    \medskip 
\vfill\eject
\noindent
{\bf The cases $s = 3$.}
     \bigskip
$$
{\beginpicture
\setcoordinatesystem units <3cm,1.2cm>
%========================================
\put{\beginpicture
\setcoordinatesystem units <.3cm,.3cm>
\put{$\ssize 1$} at -.5 3
\multiput{$\bullet$} at 0 0  1 1  2 2  3 3  2 0  4 0  6 0  3 1  5 1  4 2 /
\plot 0 0  3 3  6 0 /
\plot 1 1  2 0  4 2 /
\plot 2 2  4 0  5 1 /
\put{$\bigcirc$} at 1 1
\endpicture} at 0 0
%========================================
\put{\beginpicture
\setcoordinatesystem units <.3cm,.3cm>
\put{$\ssize 2$} at -.5 2.5
\multiput{$\bullet$} at 0 0  1 1  2 2  2 0  4 0  6 0  3 1  5 1  4 2 /
\plot 0 0  2 2 /
\plot 4 2  6 0 /
\plot 1 1  2 0  4 2 /
\plot 2 2  4 0  5 1 /
\put{$\bigcirc$} at 1 1
\put{$\bigcirc$} at 4 2
\endpicture} at 0 -2
%========================================
\put{\beginpicture
\setcoordinatesystem units <.3cm,.3cm>
\put{$\ssize 3$} at -.5 2.5
\multiput{$\bullet$} at 0 0  1 1   2 0  4 0  6 0  3 1  5 1  4 2 /
\plot 0 0  1 1 /
\plot 4 2  6 0 /
\plot 1 1  2 0  4 2 /
\plot 3 1  4 0  5 1 /
\put{$\bigcirc$} at 1 1
\put{$\bigcirc$} at 3 1
\endpicture} at 0 -4
%========================================
\put{\beginpicture
\setcoordinatesystem units <.3cm,.3cm>
\put{$\ssize 4$} at -.5 2.5
\multiput{$\bullet$} at 0 0  1 1   2 0  4 0  6 0  3 1  5 1  /
\plot 0 0  1 1  2 0  3 1  4 0  5 1  6 0  /
\put{$\bigcirc$} at 1 1
\put{$\bigcirc$} at 3 1
\put{$\bigcirc$} at 5 1
\endpicture} at 0 -6
%=====================================letzte Spalte oben
\put{\beginpicture
\setcoordinatesystem units <.8cm,.8cm>
\put{$\ssize 1$} at 0 0.5
\put{$\ssize 2$} at 0 -.5
\put{$\ssize 3$} at 1 0
\multiput{} at 1.5 0 /
\put{$\bigcirc$} at 1 0
\circulararc 300 degrees from 1.1 -.2  center at 1.5 0 
\setdots <.4mm>
\arr{0.2 0.4}{0.8 0.1}
\arr{0.2 -.4}{0.8 -.1}
\endpicture} at 3 0 
%=====================================letzt Spalte Mitte
\put{\beginpicture
\setcoordinatesystem units <1cm,1cm>
\put{$\ssize 1$} at 0.4 0
\put{$\ssize 2$} at 1.2 0
\put{$\ssize 3$} at 2 0
\multiput{$\bigcirc$} at 0.4 0  2 0 /
\arr{0.6 0}{1 0}
\circulararc 300 degrees from 2.1 -.2  center at 2.4 0 
\setdots <.4mm>
\arr{1.4 0}{1.8 0}
\endpicture} at 3 -2

%=====================================letzte Spalte unten
\put{\beginpicture
\setcoordinatesystem units <1cm,1cm>
\put{$\ssize 1$} at 0.2 0
\put{$\ssize 2$} at 1 0
\put{$\ssize 3$} at 2 0
\multiput{$\bigcirc$} at 1 0  2 0 /
\circulararc 300 degrees from 1.1 -.2  center at 1.4 0 
\circulararc 300 degrees from 2.1 -.2  center at 2.4 0 
\setdots <.4mm>
\arr{0.4 0}{0.8 0}
\endpicture} at 3 -4 
%=====================================p=2 Spalte oben
\put{\beginpicture
\setcoordinatesystem units <.8cm,.8cm>
\put{$\ssize 2$} at 0 0.5
\put{$\ssize 3$} at 0 -.5
\put{$\ssize 1$} at 1 0
\multiput{} at 1.5 0 /
\put{$\bigcirc$} at 0 -.5
\arr{0.2 -.4}{0.8 -.1}
\setdots <.4mm>
\arr{0.2 0.4}{0.8 0.1}
\circulararc 300 degrees from 1.1 -.2  center at 1.5 0 
\endpicture} at 1 0 
%=====================================p=2 Spalte Mitte
\put{\beginpicture
\setcoordinatesystem units <1cm,1cm>
\put{$\ssize 2$} at 0 0
\put{$\ssize 1$} at 1 0
\put{$\ssize 3$} at 2 0
\multiput{$\bigcirc$} at 1 0  2 0 /

\setquadratic
\plot 1.2 0.1  1.5 0.2  1.8 0.1 /
\plot 1.2 -.1  1.5 -.2  1.8 -.1 /
\arr{1.75 0.125}{1.8 0.1}
\arr{1.25 -.125}{1.2 -.1}
%\circulararc 300 degrees from 1.1 -.2  center at 1.5 0 
\setdots <.4mm>
\arr{0.2 0}{0.8 0}
\endpicture} at 1.1 -2 
%=====================================p=2.Spalte unten
\put{\beginpicture
\setcoordinatesystem units <1cm,1cm>
\put{$\ssize 2$} at 0 0
\put{$\ssize 3$} at 1 0
\put{$\ssize 1$} at 2 0
\multiput{$\bigcirc$} at 0 0  1 0 /

\arr{0.2 0}{0.8 0}
\setquadratic
\plot 1.2 0.1  1.5 0.2  1.8 0.1 /
\arr{1.75 0.125}{1.8 0.1}
%\circulararc 300 degrees from 1.1 -.2  center at 1.5 0 
\setdots <.4mm>
\plot 1.2 -.1  1.5 -.2  1.8 -.1 /
\arr{1.25 -.125}{1.2 -.1}

\endpicture} at 1.1 -4 
%=====================================p=2.Spalte oben
\put{\beginpicture
\setcoordinatesystem units <.8cm,.8cm>
\put{$\ssize 1$} at 0 0.5
\put{$\ssize 3$} at 0 -.5
\put{$\ssize 2$} at 1 0
\multiput{} at 1.5 0 /
\put{$\bigcirc$} at 0 -.5
\arr{0.2 -.4}{0.8 -.1}
\setdots <.4mm>
\arr{0.2 0.4}{0.8 0.1}
\circulararc 300 degrees from 1.1 -.2  center at 1.5 0 
\endpicture} at 2 0 
%=====================================p=2.Spalte Mitte
\put{\beginpicture
\setcoordinatesystem units <1cm,1cm>
\put{$\ssize 3$} at 0.2 0
\put{$\ssize 2$} at 1 0
\put{$\ssize 1$} at 2 0
\multiput{$\bigcirc$} at .2 0  2 0 /

\arr{0.4 0}{0.8 0}
\circulararc 300 degrees from 2.1 -.2  center at 2.4 0 
\setdots <.4mm>
\circulararc 300 degrees from 1.1 -.2  center at 1.4 0 
\endpicture} at 2 -2 
%=====================================p=2.Spalte unten
\put{\beginpicture
\setcoordinatesystem units <1cm,1cm>
\put{$\ssize 3$} at 0.4 0
\put{$\ssize 2$} at 1.2 0
\put{$\ssize 1$} at 2 0
\multiput{$\bigcirc$} at 0.4 0  1.2 0 /

\arr{0.6 0}{1 0}
\arr{1.4 0}{1.8 0}
\setdots <.4mm>
\circulararc 300 degrees from 2.1 -.2  center at 2.4 0 
\endpicture} at 2 -4 

%================================123solid
\multiput{\beginpicture
\setcoordinatesystem units <.5cm,.8cm>
\put{$\ssize 2$} at 0 0 
\put{$\ssize 1$} at 1 1
\put{$\ssize 3$} at 2 0
\multiput{$\bigcirc$} at 0 0  1 1  2 0  /
\arr{0.2 0}{1.8 0}
\arr{1.8 0.2}{1.2 0.8}
\arr{0.8 0.8}{0.2 0.2}
\endpicture} at 1.1 -6 /
%================================132solid
\multiput{\beginpicture
\setcoordinatesystem units <.5cm,.8cm>
\put{$\ssize 3$} at 0 0 
\put{$\ssize 1$} at 1 1
\put{$\ssize 2$} at 2 0
\multiput{$\bigcirc$} at 0 0  1 1  2 0  /
\arr{0.2 0}{1.8 0}
\arr{1.8 0.2}{1.2 0.8}
\arr{0.8 0.8}{0.2 0.2}
\endpicture} at 2.05 -6 /

%================================three loops solid
\multiput{\beginpicture
\setcoordinatesystem units <.5cm,.5cm>
 \put{$\ssize 1$} at 1 0
 \put{$\ssize 2$} at 2.5 0
 \put{$\ssize 3$} at 4 0
 \multiput{$\bigcirc$} at 1 0  2.5 0  4 0 /
 \circulararc 300 degrees from 1.1 -.2  center at 1.6 0 
 \circulararc 300 degrees from 2.6 -.2  center at 3.1 0 
 \circulararc 300 degrees from 4.1 -.2  center at 4.6 0 
\endpicture} at   3.05 -6 /

\plot -.4 1  3.6 1 /
\plot 0.5 1.4  0.5 -6.8 /
\setdots <1mm>
\plot .53 1  .53 -1.2  1.55  -1.2  1.55 -5  2.6  -5  2.6 -6.8 /
\setsolid
\put{$p\equiv 1\mod 3$} at 1 1.3
\put{$p\equiv 2\mod 3$} at 2 1.3
\put{$p\equiv 0\mod 3$} at 3 1.3

\put{$\ssize t=1$} at .15 -0.7
\multiput{$\ssize t=2$} at .15 -2.7  .15 -4.7   /
\put{$\ssize t=3$} at .15 -6.7

\multiput{G} at  3.35 -0.7  3.35 -2.7  3.35 -4.7  1.35 -2.7  
   1.35 -6.7  2.35 -6.7  3.35 -6.7 /

\put{$\ssize g=1,\ v=2$} at 2.95 -0.7
\put{$\ssize g=1,\ v=4$} at 2.95 -2.7
\put{$\ssize g=2,\ v=2$} at 2.95 -4.7
\put{$\ssize g=2,\ v=2$} at 1 -2.7
\multiput{$\ssize g=3,\ v=0$} at 1 -6.7  2 -6.7  3 -6.7 /

\multiput{} at  2.35 -2.7  / %the mixed cases
\put{$\ssize g=1$} at 1.9 -2.7  

\multiput{F} at 1.35 -0.7  1.35 -4.7  2.35 -0.7  2.35 -4.7  /

\endpicture}
$$
In addition of exhibiting the resolution quiver, we mention the type in question:
type G means that we deal with a Gorenstein algebra (all cycles are black); 
the type F algebras are the algebras with no black cycle in $R$, thus
the algebras of type F are $\CM$-free. There also
exist algebras which have cycles
which are black as well as cycles which are not black (thus, for $p > s$, they are
not Gorenstein algebras, 
but have non-projective Gorenstein projective modules: see the
second row with $p\equiv 2 \mod 3$).
In the left column, we mention the number $t$ of minimal projective modules,
or, equivalently, the number of indecomposable projective-injective modules. 
For the algebras which are not of type F, we add the number $g$, this is
the number of 
elementary Gorenstein projective modules, provided $p > s$. For the algebras
of type G, the number $v$ is the Gorenstein dimension (if $d$ is the  
maximal distance between a vertex and a cyclically black vertex, then $v = 2d$).
The dotted line separates the algebras $\Lambda$ with loops in $R(\Lambda)$
from those without loops (as we have mentioned in the introduction, we will show
in part II [R2] that for a connected Nakayama algebra $\Lambda$ 
there is either no loop in $R(\Lambda)$, or else all cycles in $R(\Lambda)$ are loops).
      \bigskip\bigskip
In order to arrange the possible roofs for fixed $s$, we may proceed as follows:
We consider the path algebra $\Sigma = \Sigma(s)$ of the linearly ordered quiver
of type $\Bbb A_{s+1}$, and the lattice of admissible ideals $I$ of $\Sigma$,
or, equivalently, the corresponding module categories $\mo \Sigma/I$
(using as partial ordering the inclusion functors). For any ideal $I$ we
denote by $m(I)$ the number of isomorphism classes of indecomposable
$\Sigma/I$-modules of length at least $3$. Clearly, we have
$$
 0 \le m(I) \le \binom s2.
$$
	\bigskip

Here is the {\bf roof diagram} in
the case $s = 4$ (as we will remark below, the number of admissible ideals $I$ in $\Sigma(s)$ is the Catalan number $C_s$ and $C_4 = 14$). 
    \medskip
$$
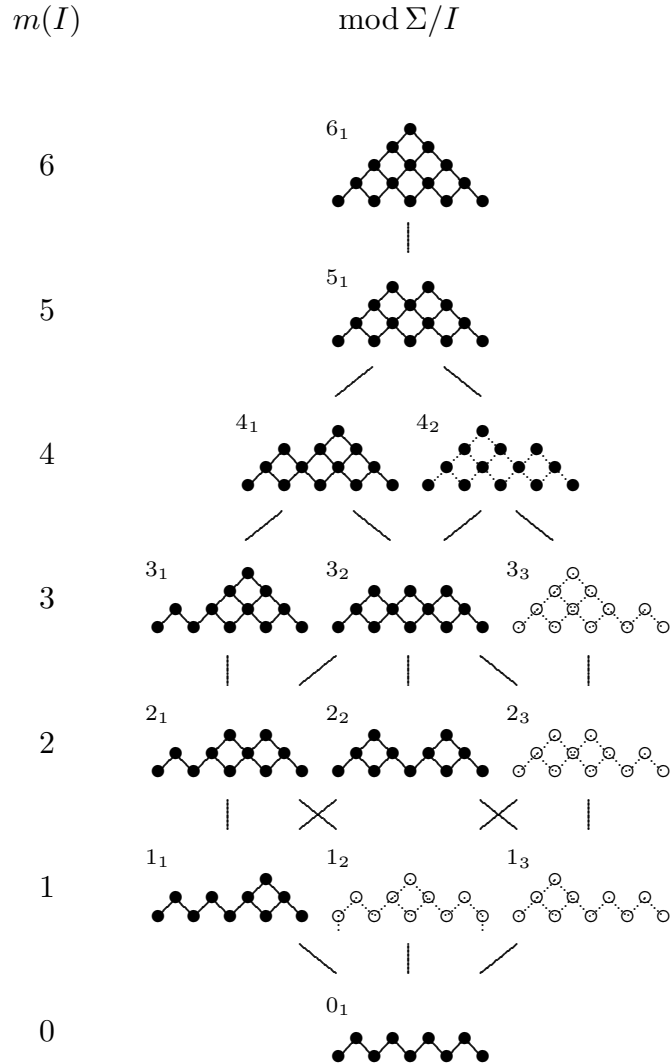
{\beginpicture
\setcoordinatesystem units <1cm,.8cm>
\put{$m(I)$} at -2 2
\put{$\mo\Sigma/I$} at 1.9 2
%==================THE ROOFS===================
\put{\beginpicture
\setcoordinatesystem units <.2cm,.2cm>
\put{$\ssize 6_1$} at 0 4
\multiput{$\bullet$} at 0 0  1 1  2 2  3 3  2 0  4 0  6 0  3 1  5 1  4 2 
 4 4  5 3  6 2  7 1  8 0 /
\plot 0 0  4 4  8 0  /
\plot 1 1  2 0  5 3 /
\plot 2 2  4 0  6 2 /
\plot 3 3  6 0  7 1 /
%\put{$\bigcirc$} at 1 1
\endpicture} at 2 0
%========================================
\put{\beginpicture
\setcoordinatesystem units <.2cm,.2cm>
\put{$\ssize 5_1$} at 0 3.5
\multiput{$\bullet$} at 0 0  1 1  2 2  3 3  2 0  4 0  6 0  3 1  5 1  4 2 
   5 3  6 2  7 1  8 0 /
\plot 0 0  3 3 /
\plot 5 3  8 0  /
\plot 1 1  2 0  5 3 /
\plot 2 2  4 0  6 2 /
\plot 3 3  6 0  7 1 /

\endpicture} at 2 -2
%========================================
\put{\beginpicture
\setcoordinatesystem units <.2cm,.2cm>
\put{$\ssize 4_1$} at 0 3.5
\multiput{$\bullet$} at 0 0  1 1  2 2   2 0  4 0  6 0  3 1  5 1  4 2 
   5 3  6 2  7 1  8 0 /
\plot 0 0  2 2 /
\plot 5 3  8 0  /
\plot 1 1  2 0  5 3 /
\plot 2 2  4 0  6 2 /
\plot 4 2  6 0  7 1 /

\endpicture} at 1 -4
%========================================
\put{\beginpicture
\setcoordinatesystem units <.2cm,.2cm>
\put{$\ssize 4_2$} at 0 3.5
\multiput{$\bullet$} at 0 0  1 1  2 2  3 3  2 0  4 0  6 0  3 1  5 1  4 2 
     6 2  7 1  8 0 /
\setdots <.5mm>
\plot 0 0  3 3 /
\plot 6 2   8 0  /
\plot 1 1  2 0  4 2 /
\plot 2 2  4 0  6 2 /
\plot 3 3  6 0  7 1 /
\endpicture} at 3 -4
%========================================
\put{\beginpicture
\setcoordinatesystem units <.2cm,.2cm>
\put{$\ssize 3_1$} at 0 3.2
\multiput{$\bullet$} at 0 0  1 1   2 0  4 0  6 0  3 1  5 1  4 2 
   5 3  6 2  7 1  8 0 /
\plot 0 0  1 1 /
\plot 5 3  8 0  /
\plot 1 1  2 0  5 3 /
\plot 3 1  4 0  6 2 /
\plot 4 2  6 0  7 1 /

\endpicture} at 0 -6
%========================================
\put{\beginpicture
\setcoordinatesystem units <.2cm,.2cm>
\put{$\ssize 3_2$} at 0 3.2
\multiput{$\bullet$} at 0 0  1 1  2 2  2 0  4 0  6 0  3 1  5 1  4 2
   6 2  7 1  8 0 /
\plot 0 0  2 2  3 1  /
\plot 6 2  8 0  /
\plot 1 1  2 0  3 1 /
\plot 3 1  4 0  6 2 /
\plot 5 1  6 0  7 1 /
\plot 3 1  4 2  5 1 /

\endpicture} at 2 -6
%========================================
\put{\beginpicture
\setcoordinatesystem units <.2cm,.2cm>
\put{$\ssize 3_3$} at 0 3.2
\multiput{$\circ$} at 0 0  1 1  2 2  3 3  2 0  4 0  6 0  3 1  5 1  4 2 
     7 1  8 0 /
\setdots <.5mm>
\plot 0 0  3 3 /
\plot 7 1   8 0  /
\plot 1 1  2 0  4 2 /
\plot 2 2  4 0  5 1 /
\plot 3 3  6 0  7 1 /
\endpicture} at 4 -6
%========================================
\put{\beginpicture
\setcoordinatesystem units <.2cm,.2cm>
\put{$\ssize 2_1$} at 0 3.2
\multiput{$\bullet$} at 0 0  1 1   2 0  4 0  6 0  3 1  5 1  4 2
   6 2  7 1  8 0 /
\plot 0 0  1 1  /
\plot 6 2  8 0  /
\plot 1 1  2 0  3 1 /
\plot 3 1  4 0  6 2 /
\plot 5 1  6 0  7 1 /
\plot 3 1  4 2  5 1 /

\endpicture} at 0 -8
%========================================
\put{\beginpicture
\setcoordinatesystem units <.2cm,.2cm>
\put{$\ssize 2_2$} at 0 3.2
\multiput{$\bullet$} at 0 0  1 1  2 2  2 0  4 0  6 0  3 1  5 1 
   6 2  7 1  8 0 /
\plot 0 0  2 2  3 1  /
\plot 6 2  8 0  /
\plot 1 1  2 0  3 1 /
\plot 3 1  4 0  6 2 /
\plot 5 1  6 0  7 1 /

\endpicture} at 2 -8
%========================================
\put{\beginpicture
\setcoordinatesystem units <.2cm,.2cm>
\put{$\ssize 2_3$} at 0 3.2
\multiput{$\circ$} at 0 0  1 1  2 2  2 0  4 0  6 0  3 1  5 1  4 2 
     7 1  8 0 /
\setdots <.5mm>
\plot 0 0  2 2 /
\plot 7 1   8 0  /
\plot 1 1  2 0  4 2 /
\plot 2 2  4 0  5 1 /
\plot 4 2  6 0  7 1 /
\endpicture} at 4 -8
%========================================
\put{\beginpicture
\setcoordinatesystem units <.2cm,.2cm>
\put{$\ssize 1_1$} at 0 3.2
\multiput{$\bullet$} at 0 0  1 1   2 0  4 0  6 0  3 1  5 1 
   6 2  7 1  8 0 /
\plot 0 0  1 1  /
\plot 6 2  8 0  /
\plot 1 1  2 0  3 1 /
\plot 3 1  4 0  6 2 /
\plot 5 1  6 0  7 1 /

\endpicture} at 0 -10
%========================================
\put{\beginpicture
\setcoordinatesystem units <.2cm,.2cm>
\put{$\ssize 1_2$} at 0 3.2
\multiput{$\circ$} at 0 0  1 1   2 0  4 0  6 0  3 1  5 1  4 2 
     7 1  8 0 /
\setdots <.5mm>
\plot 0 0  1 1 /
\plot 7 1   8 0  /
\plot 1 1  2 0  4 2 /
\plot 3 1  4 0  5 1 /
\plot 4 2  6 0  7 1 /
\setdots <.5mm>
\plot 0 0  0 -.8 /
\plot 8 0  8 -.8 /

%\put{$\bigcirc$} at 1 1
%\put{$\bigcirc$} at 5 3
\endpicture} at 2 -10
%========================================
\put{\beginpicture
\setcoordinatesystem units <.2cm,.2cm>
\put{$\ssize 1_3$} at 0 3.2
\multiput{$\circ$} at 0 0  1 1  2 2  2 0  4 0  6 0  3 1  5 1  
     7 1  8 0 /
\setdots <.5mm>
\plot 0 0  2 2 /
\plot 7 1   8 0  /
\plot 1 1  2 0  3 1 /
\plot 2 2  4 0  5 1 /
\plot 5 1   6 0  7 1 /
\endpicture} at 4 -10

%========================================
\put{\beginpicture
\setcoordinatesystem units <.2cm,.2cm>
\put{$\ssize 0_1$} at 0 2.8
\multiput{$\bullet$} at 0 0  1 1   2 0  4 0  6 0  3 1  5 1 
     7 1  8 0 /
\plot 0 0  1 1  /
\plot 1 1  2 0  3 1 /
\plot 3 1  4 0  5 1  6 0  7 1  8 0 /

\endpicture} at 2 -12

\plot 2 -0.8  2 -1.2 /
\plot 1.6 -2.8  1.2 -3.2 /
\plot 2.4 -2.8  2.8 -3.2 /
\plot .6 -4.8  .2 -5.2 /
\plot 1.4 -4.8  1.8 -5.2 /

\plot 3.2 -4.8  3.6 -5.2 /
\plot 2.8 -4.8  2.4 -5.2 /

\plot 0 -6.8  0 -7.2 /
\plot 1.2 -6.8  .8 -7.2 /
\plot 2 -6.8  2 -7.2 /
\plot 2.8 -6.8  3.2 -7.2 /
\plot 4 -6.8  4 -7.2 /

\plot 0 -8.8  0 -9.2 /
\plot 1.2 -8.8  .8 -9.2 /
\plot 2.8 -8.8  3.2 -9.2 /
\plot 4 -8.8  4 -9.2 /
\plot .8 -8.8  1.2 -9.2 /
\plot 3.2 -8.8  2.8 -9.2 /

\plot .8 -10.8  1.2 -11.2 /
\plot 2 -10.8  2 -11.2 /
\plot 3.2 -10.8  2.8 -11.2 /

\put{$6$} at -2 0
\put{$5$} at -2 -2
\put{$4$} at -2 -4
\put{$3$} at -2 -6
\put{$2$} at -2 -8
\put{$1$} at -2 -10
\put{$0$} at -2 -12

\endpicture}
$$
	\bigskip
It is sufficient to look at the roofs drawn with black bullets, the
remaining ones lead to algebras which are isomorphic to
algebras already considered: the algebras
with roof $3_3$ or $2_3$ are isomorphic to algebras with roof $3_1$
or $2_1$, respectively; the algebras with roof $1_2$ or $1_3$ are
isomorphic to algebras with roof $1_1$, always using rotations.
We may add that 
the opposite algebra of a Nakayama algebra often 
is isomorphic to the given algebra --- the only exception for $s=4$
are the cases $4_1$ and $4_2$  
(the roof $4_2$ is obtained from $4_1$ by a reflection). Of course, the opposite
algebra of a Gorenstein algebra is Gorenstein, the opposite of a CM-free algebra
is CM-free. However, the examples $4_1$ and $4_2$ show that 
the resolution quiver of the opposite
of an algebra $\Lambda$ may be quite different from $R(\Lambda)$. 
    \bigskip
{\bf Remark.} {\it For any $s$, the number of admissible ideals $I$ of $\Sigma(s)$ 
is the Catalan number $C_s = \frac 1{s+1}\binom {2s}s.$}
   \medskip
Proof: Let us rotate the Auslander-Reiten quiver of $\Sigma/I$ by
135$^o$. Then the arrows on the boundary of the Auslander-Reiten quiver
yield a monotonic path along the edges of
a grid with $(s\times s)$ square cells, starting from the lower left corner and
ending in the upper right corner (here, monotonic means that we use only
edges pointing rightwards or upwards). Here is an example with $s=3$:
$$
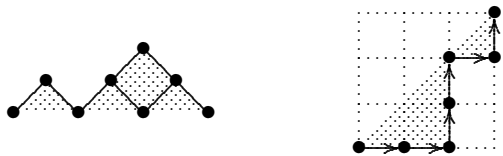
\hbox{\beginpicture
\setcoordinatesystem units <.5cm,.5cm>

%=============================================================
\put{\beginpicture
\setcoordinatesystem units <.36cm,.36cm>
\multiput{$\bullet$} at 0 0  1 1   2 0  4 0  6 0  3 1  5 1  4 2 /
\plot 0 0  1 1 /
\plot 4 2  6 0 /
\plot 1 1  2 0  4 2 /
\plot 3 1  4 0  5 1 /
%\put{$\bigcirc$} at 1 1
%\put{$\bigcirc$} at 3 1
%\setdots <.5mm>
%\plot 0 0  0 -.8 /
%\plot 6 0  6 -.8 /
\setshadegrid span <.5mm>
\vshade 0 0 0  <,z,,> 1 0 1 <z,z,,> 2 0 0 <z,z,,>  4 0 2 <z,,,>  6 0 0 /
\endpicture} at 0 0 
%=============================================================
\put{\beginpicture
\multiput{$\bullet$} at 0 0  1 0  2 0  2 1  2 2  3 2  3 3 /
\arr{0 0}{.85 0}
\arr{1 0}{1.85 0}
\arr{2 0}{2 .85}
\arr{2 1}{2 1.85}
\arr{2 2}{2.85 2}
\arr{3 2}{3 2.85}
\setdots <1mm>
\plot 0 0  0 3  3 3  3 0  0 0 /
\plot 1 0  1 3 /
\plot 2 0  2 3 /
\plot 0 1  3 1 /
\plot 0 2  3 2 /
\setshadegrid span <.5mm>
\vshade 0 0 0.05  2 0 2.05  2.01 2 2.05  3 2 3.05 /
\endpicture} at 7 0 
\endpicture} 
$$
It is well-known that the number of
such paths is just $C_s$.

     \bigskip 
Here are the Catalan numbers $C_s$ with $1\le s \le 10$.
$$
\matrix s &  1 & 2 & 3 & 4 & 5 & 6 & 7 & 8 & 9 & 10 \cr
        C_s &  1& 2& 5& 14& 42& 132& 429& 1430& 4862& 16796
\endmatrix
$$
	\bigskip
In the following table of the Nakayama algebras with $s=4$, 
the roofs are ordered slightly different:
namely, we use as first criterion the number $t$ of minimal projective
modules (or, equivalently, the number 
of indecomposable projective-injective modules). But then, the algebras
with fixed $t$ are presented in the order in which they are obtained 
in the roof diagram, reading row by row from left to right. 
   \bigskip\bigskip 
\vfill\eject
\noindent
{\bf The cases $s = 4.$}
$$
{\beginpicture
%\setcoordinatesystem units <2.6cm,1.05cm>
\setcoordinatesystem units <2.6cm,0.93cm>
%========================================
\plot -.5 1  4.7 1 /
\plot -.5 -1  4.7 -1 /
\plot -.5 -11  4.7 -11 /
\plot -.5 -17  4.7 -17 /

\plot 0.5 1.4  0.5 -18.8 /
\put{$p\equiv 1\mod 4$} at 1 1.3
\put{$p\equiv 2\mod 4$} at 2.05 1.3
\put{$p\equiv 3\mod 4$} at 3.1 1.3
\put{$p\equiv 0\mod 4$} at 4.15 1.3
\setdots <1mm>
\plot .52 1  .52 -.95  1.45  -.95  1.45 -9  2.38 -9  2.38 -16.95  3.55 -16.95  3.55 -18.8 /
\setsolid

%==================THE ROOFS===================
\put{\beginpicture
\setcoordinatesystem units <.28cm,.28cm>
\put{$\ssize 6_1$} at 0 4
\multiput{$\bullet$} at 0 0  1 1  2 2  3 3  2 0  4 0  6 0  3 1  5 1  4 2 
 4 4  5 3  6 2  7 1  8 0 /
\plot 0 0  4 4  8 0  /
\plot 1 1  2 0  5 3 /
\plot 2 2  4 0  6 2 /
\plot 3 3  6 0  7 1 /
\put{$\bigcirc$} at 1 1
\endpicture} at 0 0
%========================================
\put{\beginpicture
\setcoordinatesystem units <.28cm,.28cm>
\put{$\ssize 5_1$} at 0 3.5
\multiput{$\bullet$} at 0 0  1 1  2 2  3 3  2 0  4 0  6 0  3 1  5 1  4 2 
   5 3  6 2  7 1  8 0 /
\plot 0 0  3 3 /
\plot 5 3  8 0  /
\plot 1 1  2 0  5 3 /
\plot 2 2  4 0  6 2 /
\plot 3 3  6 0  7 1 /

\put{$\bigcirc$} at 1 1
\put{$\bigcirc$} at 5 3
\endpicture} at 0 -2
%========================================
\put{\beginpicture
\setcoordinatesystem units <.28cm,.28cm>
\put{$\ssize 4_1$} at 0 3.5
\multiput{$\bullet$} at 0 0  1 1  2 2   2 0  4 0  6 0  3 1  5 1  4 2 
   5 3  6 2  7 1  8 0 /
\plot 0 0  2 2 /
\plot 5 3  8 0  /
\plot 1 1  2 0  5 3 /
\plot 2 2  4 0  6 2 /
\plot 4 2  6 0  7 1 /
\put{$\bigcirc$} at 1 1
\put{$\bigcirc$} at 4 2

\endpicture} at 0 -4
%========================================NEU
\put{\beginpicture
\setcoordinatesystem units <.28cm,.28cm>
\put{$\ssize 4_2$} at 0 3.5
\multiput{$\bullet$} at 0 0  1 1  2 2  3 3   2 0  4 0  6 0  3 1  5 1  4 2 
     6 2  7 1  8 0 /
\plot 0 0  3 3  /
\plot 6 2   8 0  /
\plot 1 1  2 0  4 2  /
\plot 2 2  4 0  6 2 /
\plot 3 3   6 0  7 1 /
\put{$\bigcirc$} at 1 1
\put{$\bigcirc$} at 6 2

\endpicture} at 0 -6
%========================================
\put{\beginpicture
\setcoordinatesystem units <.28cm,.28cm>
\put{$\ssize 3_1$} at 0 3.5
\multiput{$\bullet$} at 0 0  1 1   2 0  4 0  6 0  3 1  5 1  4 2 
   5 3  6 2  7 1  8 0 /
\plot 0 0  1 1 /
\plot 5 3  8 0  /
\plot 1 1  2 0  5 3 /
\plot 3 1  4 0  6 2 /
\plot 4 2  6 0  7 1 /
\put{$\bigcirc$} at 1 1
\put{$\bigcirc$} at 3 1

\endpicture} at 0 -8
%========================================
\put{\beginpicture
\setcoordinatesystem units <.28cm,.28cm>
\put{$\ssize 2_2$} at 0 3.5
\multiput{$\bullet$} at 0 0  1 1  2 2  2 0  4 0  6 0  3 1  5 1 
   6 2  7 1  8 0 /
\plot 0 0  2 2  3 1  /
\plot 6 2  8 0  /
\plot 1 1  2 0  3 1 /
\plot 3 1  4 0  6 2 /
\plot 5 1  6 0  7 1 /
\put{$\bigcirc$} at 1 1
\put{$\bigcirc$} at 5 1

\endpicture} at 0 -10
%========================================
\put{\beginpicture
\setcoordinatesystem units <.28cm,.28cm>
\put{$\ssize 3_2$} at 0 3.5
\multiput{$\bullet$} at 0 0  1 1  2 2  2 0  4 0  6 0  3 1  5 1  4 2
   6 2  7 1  8 0 /
\plot 0 0  2 2  3 1  /
\plot 6 2  8 0  /
\plot 1 1  2 0  3 1 /
\plot 3 1  4 0  6 2 /
\plot 5 1  6 0  7 1 /
\plot 3 1  4 2  5 1 /
\put{$\bigcirc$} at 1 1
\put{$\bigcirc$} at 4 2
\put{$\bigcirc$} at 6 2 

\endpicture} at 0 -12
%========================================
\put{\beginpicture
\setcoordinatesystem units <.28cm,.28cm>
\put{$\ssize 2_1$} at 0 3.5
\multiput{$\bullet$} at 0 0  1 1   2 0  4 0  6 0  3 1  5 1  4 2
   6 2  7 1  8 0 /
\plot 0 0  1 1  /
\plot 6 2  8 0  /
\plot 1 1  2 0  3 1 /
\plot 3 1  4 0  6 2 /
\plot 5 1  6 0  7 1 /
\plot 3 1  4 2  5 1 /
\put{$\bigcirc$} at 1 1
\put{$\bigcirc$} at 3 1
\put{$\bigcirc$} at 6 2 

\endpicture} at 0 -14
%========================================
\put{\beginpicture
\setcoordinatesystem units <.28cm,.28cm>
\put{$\ssize 1_1$} at 0 3.5
\multiput{$\bullet$} at 0 0  1 1   2 0  4 0  6 0  3 1  5 1 
   6 2  7 1  8 0 /
\plot 0 0  1 1  /
\plot 6 2  8 0  /
\plot 1 1  2 0  3 1 /
\plot 3 1  4 0  6 2 /
\plot 5 1  6 0  7 1 /

\put{$\bigcirc$} at 1 1
\put{$\bigcirc$} at 3 1
\put{$\bigcirc$} at 5 1

\endpicture} at 0 -16
%========================================
\put{\beginpicture
\setcoordinatesystem units <.28cm,.28cm>
\put{$\ssize 0_1$} at 0 3.5
\multiput{$\bullet$} at 0 0  1 1   2 0  4 0  6 0  3 1  5 1 
    7 1  8 0 /
\plot 0 0  1 1  2 0  3 1  4 0  5 1  6 0  7 1  8 0 /
\put{$\bigcirc$} at 1 1
\put{$\bigcirc$} at 3 1
\put{$\bigcirc$} at 5 1
\put{$\bigcirc$} at 7 1

\endpicture} at 0 -17.8
%%%%%%%%%%%%%%%%%%%%%%%%%%%%%%%%%%%%%%%%%%%%%%%%%%
%=====================================1 a=1
\put{\beginpicture
\setcoordinatesystem units <.8cm,.8cm>
\put{$\ssize 2$} at 0 0.5
\put{$\ssize 3$} at 0 0
\put{$\ssize 4$} at 0 -.5
\put{$\ssize 1$} at 1 0
\multiput{} at 1.5 0 /
\put{$\bigcirc$} at 0 -.5
\arr{0.2 -.4}{0.8 -.1}
\setdots <.4mm>
\circulararc 300 degrees from 1.1 -.2  center at 1.5 0 
\arr{0.2 0.4}{0.8 0.1}
\arr{0.2 0}{0.8 0}
\endpicture} at 1 0 
%=====================================1 a=2
\put{\beginpicture
\setcoordinatesystem units <.8cm,.8cm>
\put{$\ssize 1$} at 0 0.5
\put{$\ssize 3$} at 0 0
\put{$\ssize 4$} at 0 -.5
\put{$\ssize 2$} at 1 0
\multiput{} at 1.5 0 /
\put{$\bigcirc$} at 0 -.5
\arr{0.2 -.4}{0.8 -.1}
\setdots <.4mm>
\circulararc 300 degrees from 1.1 -.2  center at 1.5 0 
\arr{0.2 0.4}{0.8 0.1}
\arr{0.2 0}{0.8 0}
\endpicture} at 2 0 
%=====================================1 a=3
\put{\beginpicture
\setcoordinatesystem units <.8cm,.8cm>
\put{$\ssize 1$} at 0 0.5
\put{$\ssize 2$} at 0 0
\put{$\ssize 4$} at 0 -.5
\put{$\ssize 3$} at 1 0
\multiput{} at 1.5 0 /
\put{$\bigcirc$} at 0 -.5
\arr{0.2 -.4}{0.8 -.1}
\setdots <.4mm>
\circulararc 300 degrees from 1.1 -.2  center at 1.5 0 
\arr{0.2 0.4}{0.8 0.1}
\arr{0.2 0}{0.8 0}
\endpicture} at 3 0 

%=====================================1, a=0
\put{\beginpicture
\setcoordinatesystem units <.8cm,.8cm>
\put{$\ssize 1$} at 0 0.5
\put{$\ssize 2$} at 0 0
\put{$\ssize 3$} at 0 -.5
\put{$\ssize 4$} at 1 0
\multiput{} at 1.5 0 /
\put{$\bigcirc$} at 1 0
\circulararc 300 degrees from 1.1 -.2  center at 1.5 0 
\setdots <.4mm>
\arr{0.2 0.4}{0.8 0.1}
\arr{0.2 0}{0.8 0}
\arr{0.2 -.4}{0.8 -.1}
\endpicture} at 4.15 0

%%%%%%%%%%%%%%%%%%%%%%%%%%%%%%%%%%%%

%=====================================2 a=1
\put{\beginpicture
\setcoordinatesystem units <.8cm,.8cm>
\put{$\ssize 2$} at 0 .5
\put{$\ssize 3$} at 0 -.5
\put{$\ssize 1$} at 1 0
\put{$\ssize 4$} at 2 0
\multiput{$\bigcirc$} at 1 0  2 0 /

\setquadratic
\plot 1.2 0.1  1.5 0.2  1.8 0.1 /
\plot 1.2 -.1  1.5 -.2  1.8 -.1 /
\arr{1.75 0.125}{1.8 0.1}
\arr{1.25 -.125}{1.2 -.1}
%\circulararc 300 degrees from 1.1 -.2  center at 1.5 0 
\setdots <.4mm>
\arr{0.2 0.4}{0.8 0.1}
\arr{0.2 -.4}{0.8 -.1}
\endpicture} at 1 -2 

%=====================================2 a=2
\put{\beginpicture
\setcoordinatesystem units <.8cm,.8cm>
\put{$\ssize 3$} at 0 0.5
\put{$\ssize 1$} at 2.1 0
\put{$\ssize 4$} at 0 -.5
\put{$\ssize 2$} at 1 0
\multiput{} at 1.5 0 /
\multiput{$\bigcirc$} at 0 -.5  2.1 0 /
\arr{0.2 -.4}{0.8 -.1}
\circulararc 300 degrees from 2.2 -.2  center at 2.5 0 
\setdots <.4mm>
\circulararc 300 degrees from 1.1 -.2  center at 1.4 0 
\arr{0.2 0.4}{0.8 0.1}
%\arr{0.2 0}{0.8 0}
\endpicture} at 2 -2 

%=====================================2 a=3
\put{\beginpicture
\setcoordinatesystem units <.8cm,.8cm>
\put{$\ssize 2$} at 0 0.5
\put{$\ssize 1$} at -0.8 0.5
\put{$\ssize 4$} at 0 -.5
\put{$\ssize 3$} at 1 0
\multiput{} at 1.5 0 /
\multiput{$\bigcirc$} at 0 -.5  -0.8 0.5 /
\arr{0.2 -.4}{0.8 -.1}
\arr{-.6 0.5}{-.2 0.5}
\setdots <.4mm>
\circulararc 300 degrees from 1.1 -.2  center at 1.5 0 
\arr{0.2 0.4}{0.8 0.1}
\endpicture} at 3.1 -2 

%=====================================2 a=0
\put{\beginpicture
\setcoordinatesystem units <.8cm,.8cm>
\put{$\ssize 3$} at 0 0.5
\put{$\ssize 1$} at -0.8 0.5
\put{$\ssize 2$} at 0 -.5
\put{$\ssize 4$} at 1 0
\multiput{} at 1.5 0 /
\multiput{$\bigcirc$} at 1 0  -0.8 0.5 /
\arr{-.6 0.5}{-.2 0.5}
\circulararc 300 degrees from 1.1 -.2  center at 1.5 0 
\setdots <.4mm>
\arr{0.2 -.4}{0.8 -.1}
\arr{0.2 0.4}{0.8 0.1}
\endpicture} at 4.15 -2

%%%%%%%%%%%%%%%%%%%%%%%%%%%%%%%%%%
%=====================================3 a=1
\put{\beginpicture
\setcoordinatesystem units <.8cm,.8cm>
\put{} at -1 0
\put{$\ssize 2$} at 0 0
\put{$\ssize 3$} at 0 1
\put{$\ssize 4$} at 1 0 
\put{$\ssize 1$} at 1 1
\multiput{$\bigcirc$} at 1 0  0 0 /
\arr{0.89 0.22}{0.9 0.2}
\arr{1.11 0.78}{1.1 0.8}
\setquadratic
\plot 1.1 0.2 1.2 0.5 1.1 0.8 /
\arr{0.2 0}{0.8 0}
\setdots <.4mm>
\plot 0.89 0.2 0.8 0.5 0.9 0.8 /
\arr{0.2 1}{0.8 1}
\endpicture} at 1 -4

%=====================================3 a=2
\put{\beginpicture
\setcoordinatesystem units <.8cm,.8cm>
\put{$\ssize 3$} at 0.4 0.5
\put{$\ssize 4$} at 0.4 -.5
\put{$\ssize 2$} at 1.2 0
\put{$\ssize 1$} at 2 0
\multiput{$\bigcirc$} at 0.4 -.5  1.2 0 /
\arr{0.6 -.4}{1 -.1}
\arr{1.4 0}{1.8 0}
\setdots <.4mm>
\circulararc 300 degrees from 2.1 -.2  center at 2.5 0 
\arr{0.6 0.4}{1 0.1}
\endpicture} at 1.9 -4

%=====================================3 a=3
\put{\beginpicture
\setcoordinatesystem units <.8cm,.8cm>
\put{$\ssize 1$} at 0.2 0
\put{$\ssize 2$} at 1 0
\put{$\ssize 4$} at 2.2 0
\put{$\ssize 3$} at 3 0
%\multiput{} at 1.5 0 /
\multiput{$\bigcirc$} at 1 0  2.2 0 /
\arr{2.4 0}{2.8 0}
\circulararc 300 degrees from 1.1 -.2  center at 1.4 0 
\setdots <.4mm>
\circulararc 300 degrees from 3.1 -.2  center at 3.4 0 
\arr{0.4 0}{0.8 0}
%\arr{0.2 0}{0.8 0}
\endpicture} at 3 -4

%=====================================3 a=0
\put{\beginpicture
\setcoordinatesystem units <.8cm,.8cm>
\put{$\ssize 1$} at 0.4 0.5
\put{$\ssize 2$} at 0.4 -.5
\put{$\ssize 3$} at 1.2 0
\put{$\ssize 4$} at 2 0
\multiput{$\bigcirc$} at 0.4 -.5  2 0 /
\arr{0.6 -.4}{1 -.1}
\circulararc 300 degrees from 2.1 -.2  center at 2.5 0 
\setdots <.4mm>
\arr{0.6 0.4}{1 0.1}
\arr{1.4 0}{1.8 0}
\endpicture} at 4.15 -4 

%%%%%%%%%%%%%%%%%%%%%%%%%%%%%%%%%%%%%%%INSERTED
%=====================================NEU
\put{\beginpicture
\setcoordinatesystem units <.8cm,.8cm>
\put{$\ssize 2$} at 0 .5
\put{$\ssize 4$} at 0 -.5
\put{$\ssize 1$} at 1 0
\put{$\ssize 3$} at 2 0
\multiput{$\bigcirc$} at 1 0  0 -.5 /

\setquadratic
\plot 1.2 0.1  1.5 0.2  1.8 0.1 /
\arr{1.75 0.125}{1.8 0.1}
\arr{1.25 -.125}{1.2 -.1}
\arr{0.2 -.4}{0.8 -.1}
\setdots <.4mm>
\arr{0.2 0.4}{0.8 0.1}
\plot 1.2 -.1  1.5 -.2  1.8 -.1 /
\endpicture} at 1 -6 

%=====================================NEU
\put{\beginpicture
\setcoordinatesystem units <.8cm,.8cm>
\put{$\ssize 4$} at 0 0.5
\put{$\ssize 1$} at -0.8 0.5
\put{$\ssize 3$} at 0 -.5
\put{$\ssize 2$} at 1 0
\multiput{} at 1.5 0 /
\multiput{$\bigcirc$} at 0 0.5  -0.8 0.5 /
\arr{-.6 0.5}{-.2 0.5}
\arr{0.2 0.4}{0.8 0.1}
\setdots <.4mm>
\circulararc 300 degrees from 1.1 -.2  center at 1.5 0 
\arr{0.2 -.4}{0.8 -.1}
\endpicture} at 1.93 -6

%====================================NEU
\put{\beginpicture
\setcoordinatesystem units <.8cm,.8cm>
\put{$\ssize 2$} at 0 0.5
\put{$\ssize 1$} at 2.1 0
\put{$\ssize 4$} at 0 -.5
\put{$\ssize 3$} at 1 0
\multiput{} at 1.5 0 /
\multiput{$\bigcirc$} at 0 -.5  2.1 0 /
\arr{0.2 -.4}{0.8 -.1}
\circulararc 300 degrees from 2.2 -.2  center at 2.5 0 
\setdots <.4mm>
\circulararc 300 degrees from 1.1 -.2  center at 1.4 0 
\arr{0.2 0.4}{0.8 0.1}
%\arr{0.2 0}{0.8 0}
\endpicture} at 3 -6 

%=====================================NEU
\put{\beginpicture
\setcoordinatesystem units <.8cm,.8cm>
\put{$\ssize 2$} at 0 0.5
\put{$\ssize 1$} at -0.8 0.5
\put{$\ssize 3$} at 0 -.5
\put{$\ssize 4$} at 1 0
\multiput{} at 1.5 0 /
\multiput{$\bigcirc$} at 1 0  -0.8 0.5 /
\arr{-.6 0.5}{-.2 0.5}
\circulararc 300 degrees from 1.1 -.2  center at 1.5 0 
\setdots <.4mm>
\arr{0.2 -.4}{0.8 -.1}
\arr{0.2 0.4}{0.8 0.1}
\endpicture} at 4.15 -6

%%%%%%%%%%%%%%%%%%%%%%%%%%%%%%%%%%%%%%%
%=====================================4 a=1
\put{\beginpicture
\setcoordinatesystem units <.8cm,.8cm>
\put{$\ssize 2$} at 0 .5
\put{$\ssize 3$} at 0 -.5
\put{$\ssize 4$} at 1 0
\put{$\ssize 1$} at 2 0
\multiput{$\bigcirc$} at 1 0  0 -.5 /

\setquadratic
\plot 1.2 0.1  1.5 0.2  1.8 0.1 /
\arr{1.75 0.125}{1.8 0.1}
\arr{1.25 -.125}{1.2 -.1}
\arr{0.2 -.4}{0.8 -.1}
\setdots <.4mm>
\arr{0.2 0.4}{0.8 0.1}
\plot 1.2 -.1  1.5 -.2  1.8 -.1 /
\endpicture} at 1 -8 

%=====================================4 a=2
\put{\beginpicture
\setcoordinatesystem units <.8cm,.8cm>
\put{$\ssize 2$} at 0 0.5
\put{$\ssize 4$} at -0.8 0.5
\put{$\ssize 3$} at 0 -.5
\put{$\ssize 1$} at 1 0
\multiput{} at 1.5 0 /
\multiput{$\bigcirc$} at 0 -.5  -0.8 0.5 /
\arr{0.2 -.4}{0.8 -.1}
\arr{-.6 0.5}{-.2 0.5}
\setdots <.4mm>
\circulararc 300 degrees from 1.1 -.2  center at 1.5 0 
\arr{0.2 0.4}{0.8 0.1}

\endpicture} at 1.9 -8 

%=====================================4 a=3
\put{\beginpicture
\setcoordinatesystem units <.8cm,.8cm>
\put{$\ssize 3$} at 0 0.5
\put{$\ssize 4$} at -0.8 0.5
\put{$\ssize 1$} at 0 -.5
\put{$\ssize 2$} at 1 0
\multiput{} at 1.5 0 /
\multiput{$\bigcirc$} at 0 0.5  -0.8 0.5 /
\arr{-.6 0.5}{-.2 0.5}
\arr{0.2 0.4}{0.8 0.1}
\setdots <.4mm>
\circulararc 300 degrees from 1.1 -.2  center at 1.5 0 
\arr{0.2 -.4}{0.8 -.1}
\endpicture} at 3 -8

%=====================================4 a=0
\put{\beginpicture
\setcoordinatesystem units <.8cm,.8cm>
\put{$\ssize 1$} at 0 0.5
\put{$\ssize 4$} at 2.1 0
\put{$\ssize 2$} at 0 -.5
\put{$\ssize 3$} at 1 0
\multiput{} at 1.5 0 /
\multiput{$\bigcirc$} at 1 0  2.1 0 /
\circulararc 300 degrees from 2.2 -.2  center at 2.5 0 
\circulararc 300 degrees from 1.1 -.2  center at 1.4 0 
\setdots <.4mm>
\arr{0.2 -.4}{0.8 -.1}
\arr{0.2 0.4}{0.8 0.1}
%\arr{0.2 0}{0.8 0}
\endpicture} at 4.15 -8 

%%%%%%%%%%%%%%%%%%%%%%%%%%%%%%%%%%%%%%%%%%%%%%%%%%

%=====================================5 a=1
\put{\beginpicture
\setcoordinatesystem units <.8cm,.8cm>
\put{} at -1 0
\put{$\ssize 2$} at 0 0
\put{$\ssize 4$} at 0 1
\put{$\ssize 3$} at 1 0 
\put{$\ssize 1$} at 1 1
\multiput{$\bigcirc$} at 0 0  0 1 /
\arr{0.2 0}{0.8 0}
\arr{0.2 1}{0.8 1}
\arr{0.89 0.22}{0.9 0.2}
\arr{1.11 0.78}{1.1 0.8}
\setdots <.4mm>
\setquadratic
\plot 0.89 0.2 0.8 0.5 0.9 0.8 /
\plot 1.1 0.2 1.2 0.5 1.1 0.8 /
\endpicture} at 1 -10 

%=====================================5 a=2
\put{\beginpicture
\setcoordinatesystem units <.8cm,.8cm>
\put{} at -1 0
\put{$\ssize 3$} at 0 0
\put{$\ssize 1$} at 0 1
\put{$\ssize 2$} at 1 0 
\put{$\ssize 4$} at 1 1
\multiput{$\bigcirc$} at 1 0  1 1 /
\arr{0.89 0.22}{0.9 0.2}
\arr{1.11 0.78}{1.1 0.8}
\setquadratic
\plot 0.89 0.2 0.8 0.5 0.9 0.8 /
\plot 1.1 0.2 1.2 0.5 1.1 0.8 /
\setdots <.4mm>
\arr{0.2 0}{0.8 0}
\arr{0.2 1}{0.8 1}
\endpicture} at 1.9 -10

%=====================================5 a=3
\put{\beginpicture
\setcoordinatesystem units <.8cm,.8cm>
\put{$\ssize 2$} at 0.2 0
\put{$\ssize 1$} at 1 0
\put{$\ssize 4$} at 2.2 0
\put{$\ssize 3$} at 3 0
%\multiput{} at 1.5 0 /
\multiput{$\bigcirc$} at .2 0  2.2 0 /
\arr{2.5 0}{2.8 0}
\arr{0.5 0}{.8 0}
\setdots <.4mm>
\circulararc 300 degrees from 1.1 -.2  center at 1.4 0 
\circulararc 300 degrees from 3.1 -.2  center at 3.4 0 
%\arr{0.2 0}{0.8 0}
\endpicture} at 2.9 -10

%=====================================5 a=0
\put{\beginpicture
\setcoordinatesystem units <.8cm,.8cm>
\put{$\ssize 1$} at 0.4 0
\put{$\ssize 2$} at 1.2 0
\put{$\ssize 3$} at 2.2 0
\put{$\ssize 4$} at 3 0
%\multiput{} at 1.5 0 /
\multiput{$\bigcirc$} at 1.2 0  3 0 /
\circulararc 300 degrees from 1.3 -.2  center at 1.6 0 
\circulararc 300 degrees from 3.1 -.2  center at 3.4 0 
\setdots <.4mm>
\arr{2.4 0}{2.8 0}
\arr{0.6 0}{1 0}
%\arr{0.2 0}{0.8 0}
\endpicture} at 4.12 -10

%%%%%%%%%%%%%%%%%%%%%%%%%%%%%%%%%%%%%%%%%%%%%%%%%%%%%%%%

%=====================================6 a=1
\put{\beginpicture
\setcoordinatesystem units <.8cm,.8cm>
\put{$\ssize 2$} at -.6 0
\put{$\ssize 4$} at 0.2 0
\put{$\ssize 1$} at 1 0
\put{$\ssize 3$} at 2 0
\multiput{$\bigcirc$} at -.6 0  0.2 0  1 0  /

\setquadratic
\plot 1.2 0.1  1.5 0.2  1.8 0.1 /
\arr{1.75 0.125}{1.8 0.1}
\arr{1.25 -.125}{1.2 -.1}
\arr{0.4 0}{0.8 0}
\arr{-.4 0}{0 0}
\setdots <.4mm>
\plot 1.2 -.1  1.5 -.2  1.8 -.1 /
\endpicture} at 1.05 -12

%=====================================6 a=2
\put{\beginpicture
\setcoordinatesystem units <.8cm,.8cm>
\put{} at -1 0
\put{$\ssize 3$} at 0.2 0
\put{$\ssize 2$} at 1 0
\put{$\ssize 1$} at 2 0.5
\put{$\ssize 4$} at 2 -.5
\multiput{$\bigcirc$} at 1 0  2 0.5  2 -.5  /

\arr{1.2 0.1}{1.8 0.4}
\arr{2 0.3}{2 -.3}
\arr{1.8 -.4}{1.2 -.1}
\setdots <.4mm>
\arr{0.4 0}{0.8 0}

\endpicture} at 1.82 -12

%=====================================6 a=3
\put{\beginpicture
\setcoordinatesystem units <.8cm,.8cm>
\put{$\ssize 1$} at 0 0
\put{$\ssize 2$} at 1.1 0
\put{$\ssize 4$} at 2.1 0
\put{$\ssize 3$} at 2.9 0
\multiput{$\bigcirc$} at  0 0  1.1 0  2.1 0 /
\circulararc 300 degrees from 0.1 -.2  center at 0.4 0 
\circulararc 300 degrees from 1.2 -.2  center at 1.5 0 
\arr{2.3 0}{2.7 0}
\setdots <.4mm>
\circulararc 300 degrees from 3 -.2  center at 3.3 0 
\endpicture} at 2.9 -12

%=====================================6 a=0
\put{\beginpicture
\setcoordinatesystem units <.8cm,.8cm>
\put{$\ssize 1$} at 0.6 0
\put{$\ssize 2$} at 1.4 0
\put{$\ssize 3$} at 2.2 0
\put{$\ssize 4$} at 3 0
\multiput{$\bigcirc$} at  0.6 0  1.4 0  3 0 /
\arr{1.6 0}{2 0}
\arr{0.8 0}{1.2 0}
\circulararc 300 degrees from 3.1 -.2  center at 3.4 0 
\setdots <.4mm>
\arr{2.4 0}{2.8 0}
\endpicture} at 4.12 -12
%%%%%%%%%%%%%%%%%%%%%%%%%%%%%%%%%%%%%%%%%%%%
%=====================================7 a=1
\put{\beginpicture
\setcoordinatesystem units <.8cm,.8cm>
\put{} at -1 0
\put{$\ssize 2$} at 0.2 0
\put{$\ssize 4$} at 1 0
\put{$\ssize 1$} at 2 0.5
\put{$\ssize 3$} at 2 -.5
\multiput{$\bigcirc$} at 1 0  2 0.5  2 -.5  /

\arr{1.2 0.1}{1.8 0.4}
\arr{2 0.3}{2 -.3}
\arr{1.8 -.4}{1.2 -.1}
\setdots <.4mm>
\arr{0.4 0}{0.8 0}

\endpicture} at .86 -14

%=====================================7 a=2
\put{\beginpicture
\setcoordinatesystem units <.8cm,.8cm>
\put{} at -1 0
\put{$\ssize 3$} at 0.2 0
\put{$\ssize 1$} at 1 0
\put{$\ssize 4$} at 2 0.5
\put{$\ssize 2$} at 2 -.5
\multiput{$\bigcirc$} at 0.2 0  1 0  2 0.5  /

\arr{1.2 0.1}{1.8 0.4}
\arr{2 0.3}{2 -.3}
\arr{0.4 0}{0.8 0}
\setdots <.4mm>
\arr{1.8 -.4}{1.2 -.1}

\endpicture} at 1.84 -14

%=====================================7 a=3
\put{\beginpicture
\setcoordinatesystem units <.8cm,.8cm>
\put{$\ssize 4$} at 0.7 0
\put{$\ssize 3$} at 1.4 0
\put{$\ssize 2$} at 2.2 0
\put{$\ssize 1$} at -.4 0
\multiput{$\bigcirc$} at  0.7 0  1.4 0  -.4 0 /
\arr{0.9 0}{1.2 0}
\arr{1.6 0}{2 0}
\circulararc 300 degrees from -.4 -.2  center at -.1 0 
\setdots <.4mm>
\circulararc 300 degrees from 2.3 -.2  center at 2.6 0 
\endpicture} at 2.9 -14

%=====================================7 a=0
\put{\beginpicture
\setcoordinatesystem units <.8cm,.8cm>
\put{$\ssize 1$} at 0.6 0
\put{$\ssize 2$} at 1.4 0
\put{$\ssize 3$} at 2.2 0
\put{$\ssize 4$} at 3.3 0
\multiput{$\bigcirc$} at  0.6 0  2.2 0  3.3 0 /
\arr{0.8 0}{1.2 0}
\circulararc 300 degrees from 2.3 -.2  center at 2.6 0 
\circulararc 300 degrees from 3.4 -.2  center at 3.7 0 
\setdots <.4mm>
\arr{1.6 0}{2 0}
\endpicture} at 4.1 -14

%%%%%%%%%%%%%%%%%%%%%%%%%%%%%%%%%%%%%%%%%%%%%%%%%%%%%%
%=====================================8 a=1
\put{\beginpicture
\setcoordinatesystem units <.8cm,.8cm>
\put{} at -1 0
\put{$\ssize 2$} at 0.2 0
\put{$\ssize 3$} at 1 0
\put{$\ssize 4$} at 2 0.5
\put{$\ssize 1$} at 2 -.5
\multiput{$\bigcirc$} at 0.2 0 1 0  2 0.5  /

\arr{1.2 0.1}{1.8 0.4}
\arr{2 0.3}{2 -.3}
\arr{0.4 0}{0.8 0}
\setdots <.4mm>
\arr{1.8 -.4}{1.2 -.1}

\endpicture} at 0.86 -16 

%=====================================8 a=2
\put{\beginpicture
\setcoordinatesystem units <.8cm,.8cm>
\put{$\ssize 3$} at -.6 0
\put{$\ssize 1$} at 0.2 0
\put{$\ssize 4$} at 1 0
\put{$\ssize 2$} at 2 0
\multiput{$\bigcirc$} at -.6 0  1 0  2 0  /

\setquadratic
\plot 1.2 0.1  1.5 0.2  1.8 0.1 /
\arr{1.75 0.125}{1.8 0.1}
\arr{1.25 -.125}{1.2 -.1}
\arr{-.4 0}{0 0}
\plot 1.2 -.1  1.5 -.2  1.8 -.1 /
\setdots <.4mm>
\arr{0.4 0}{0.8 0}
\endpicture} at 1.88 -16

%=====================================8 a=3
\put{\beginpicture
\setcoordinatesystem units <.8cm,.8cm>
\put{$\ssize 4$} at 0.6 0
\put{$\ssize 3$} at 1.4 0
\put{$\ssize 2$} at 2.2 0
\put{$\ssize 1$} at 3 0
\multiput{$\bigcirc$} at  0.6 0  1.4 0  2.2 0 /
\arr{1.6 0}{2 0}
\arr{0.8 0}{1.2 0}
\arr{2.4 0}{2.8 0}
\setdots <.4mm>
\circulararc 300 degrees from 3.1 -.2  center at 3.4 0 
\endpicture} at 2.9 -16
%=====================================8 a=0
\put{\beginpicture
\setcoordinatesystem units <.8cm,.8cm>
\put{$\ssize 1$} at 0.3 0
\put{$\ssize 2$} at 1.1 0
\put{$\ssize 3$} at 2.2 0
\put{$\ssize 4$} at 3.3 0
\multiput{$\bigcirc$} at  1.1 0  2.2 0  3.3 0 /
\circulararc 300 degrees from 1.2 -.2  center at 1.5 0 
\circulararc 300 degrees from 2.3 -.2  center at 2.6 0 
\circulararc 300 degrees from 3.4 -.2  center at 3.7 0 
\setdots <.4mm>
\arr{0.5 0}{0.9 0}
\endpicture} at 4.1 -16

%%%%%%%%%%%%%%%%%%%%%%%%%%%%%%%%%%%%%%%%%%%%%%%%%%%%%%%%%%%

\put{$\ssize t=1$} at .15 -0.75
\multiput{$\ssize t=2$} at .15 -2.75  .15 -4.75  .15 -8.75  .15 -10.75 /
\multiput{$\ssize t=3$} at .15 -12.75  .15 -14.75  .15 -16.75 /

\multiput{G} at  4.5 -0.7  4.5 -2.7  
 4.5 -4.8
 4.5 -6.8
 4.5 -8.7 
 4.5 -10.7 
 4.5 -12.7 
 4.5 -14.7 
 4.5 -16.7 
 1.3 -2.7  
 2.3 -10.7 
 2.3 -12.7 
 1.3 -14.7 
 2.3 -16.7 
 1.3 -18.7
 2.3 -18.7
 3.3 -18.7 
 4.5 -18.7
/
\put{$\ssize g=1,\ v=2$} at 4.05 -0.7
\put{$\ssize g=1,\ v=4$} at 4.05 -2.7
\put{$\ssize g=1,\ v=4$} at 4.05 -4.8
\put{$\ssize g=1,\ v=4$} at 4.05 -6.8 %NEU
\put{$\ssize g=2,\ v=2$} at 4.05 -8.7
\put{$\ssize g=2,\ v=2$} at 4.05 -10.75
\put{$\ssize g=1,\ v=6$} at 4.05 -12.7
\put{$\ssize g=2,\ v=4$} at 4.05 -14.7
\put{$\ssize g=3,\ v=2$} at 4.05 -16.7

\put{$\ssize g=2,\ v=2$} at  0.9 -2.7  
\put{$\ssize g=2,\ v=2$} at  1.9 -10.75 
\put{$\ssize g=3,\ v=2$} at  1.9 -12.7 
\put{$\ssize g=3,\ v=2$} at  0.9 -14.7 
\put{$\ssize g=2,\ v=4$} at  1.9 -16.7 

\multiput{$\ssize g=4,\ v=0$} at  0.9 -16.7  1.9 -18.7  2.9 -18.7  4.05 -18.7 /

\multiput{} at  2.3 -2.7  %the mixed cases
  3.35 -4.7 
  3.35 -12.7 
  3.35 -14.7 
 /
\put{$\ssize g=1$} at  2 -2.7  
\put{$\ssize g=1$} at  3 -4.8  
\put{$\ssize g=1$} at  3 -6.8  %NEU
\put{$\ssize g=2$} at  3 -12.7 
\put{$\ssize g=1$} at  3 -14.7

\multiput{F} at 1.25 -0.7 
  1.25 -4.8 
  1.25 -6.8   %NEU
  1.25 -8.7 
  1.25 -10.7 
  1.25 -12.7 
  1.25 -16.7 
  2.3 -0.7  
  2.3 -4.8  
  2.3 -6.8  %NEU
  2.3 -8.7  
  2.3 -14.7  
  3.35 -0.7  
  3.35 -2.7  
  3.35 -8.7  
  3.35 -10.7  
  3.35 -16.7  /

%================================loops
\multiput{\beginpicture
\setcoordinatesystem units <.6cm,.6cm>
 \put{$\ssize 1$} at 1 0
 \put{$\ssize 2$} at 2.3 0
 \put{$\ssize 3$} at 3.6 0
 \put{$\ssize 4$} at 4.9 0
 \multiput{} at 1.5 0 /
 \multiput{$\bigcirc$} at 1 0 2.3 0  3.6 0  4.9 0 /
 \circulararc 300 degrees from 1.1 -.2  center at 1.5 0 
 \circulararc 300 degrees from 2.4 -.2  center at 2.8 0 
 \circulararc 300 degrees from 3.7 -.2  center at 4.1 0 
 \circulararc 300 degrees from 5.0 -.2  center at 5.4 0 
\endpicture} at 4.1 -17.8 /

%================================1234
\multiput{\beginpicture
\setcoordinatesystem units <.8cm,.8cm>
\put{$\ssize 1$} at 0 1
\put{$\ssize 2$} at 0 0 
\put{$\ssize 3$} at 1 0
\put{$\ssize 4$} at 1 1
\multiput{$\bigcirc$} at 0 0  1 0  1 1   0 1 /
\setsolid
\arr{0 0.8}{0 0.2}
\arr{0.2 0}{0.8 0}
\arr{0.8 1}{0.2 1}
\arr{1 0.2}{1 0.8}
\endpicture} at 1 -17.8 /
%================================1432
\multiput{\beginpicture
\setcoordinatesystem units <.8cm,.8cm>
\put{$\ssize 1$} at 0 1
\put{$\ssize 4$} at 0 0 
\put{$\ssize 2$} at 1 0
\put{$\ssize 4$} at 1 1
\multiput{$\bigcirc$} at 0 0  1 0  1 1   0 1 /
\setsolid
\arr{0 0.8}{0 0.2}
\arr{0.2 0}{0.8 0}
\arr{0.8 1}{0.2 1}
\arr{1 0.2}{1 0.8}
\endpicture} at 3 -17.8 /

%================================14solid
\multiput{\beginpicture
\setcoordinatesystem units <.8cm,.8cm>
\put{$\ssize 3$} at 1 0 
\put{$\ssize 1$} at 1 1
\multiput{$\bigcirc$} at 1 0  1 1 /
\arr{0.89 0.22}{0.9 0.2}
\arr{1.11 0.78}{1.1 0.8}
\setquadratic
\plot 1.1 0.2 1.2 0.5 1.1 0.8 /
\plot 0.89 0.2 0.8 0.5 0.9 0.8 /

\endpicture} at 1.9 -17.8 /

%================================14solid
\multiput{\beginpicture
\setcoordinatesystem units <.8cm,.8cm>
\put{$\ssize 4$} at 1 0 
\put{$\ssize 2$} at 1 1
\multiput{$\bigcirc$} at 1 0  1 1 /
\arr{0.89 0.22}{0.9 0.2}
\arr{1.11 0.78}{1.1 0.8}
\setquadratic
\plot 1.1 0.2 1.2 0.5 1.1 0.8 /
\plot 0.89 0.2 0.8 0.5 0.9 0.8 /

\endpicture} at 2.2 -17.8 /

\endpicture}
$$

%%%%%%%%%%%%%%s=5
$$
{\beginpicture
\setcoordinatesystem units <2.4cm,.86cm>
%========================================
\put{{\bf The cases $s = 5$:}  the cycles of the resolution quivers.} [l] at -.5 2 
\plot -.5 1  5.4 1 /
\plot -.5 -1  5.4 -1 /
\plot -.5 -13  5.4 -13 /
\plot 0.55 1.4  0.55 -18.7 /

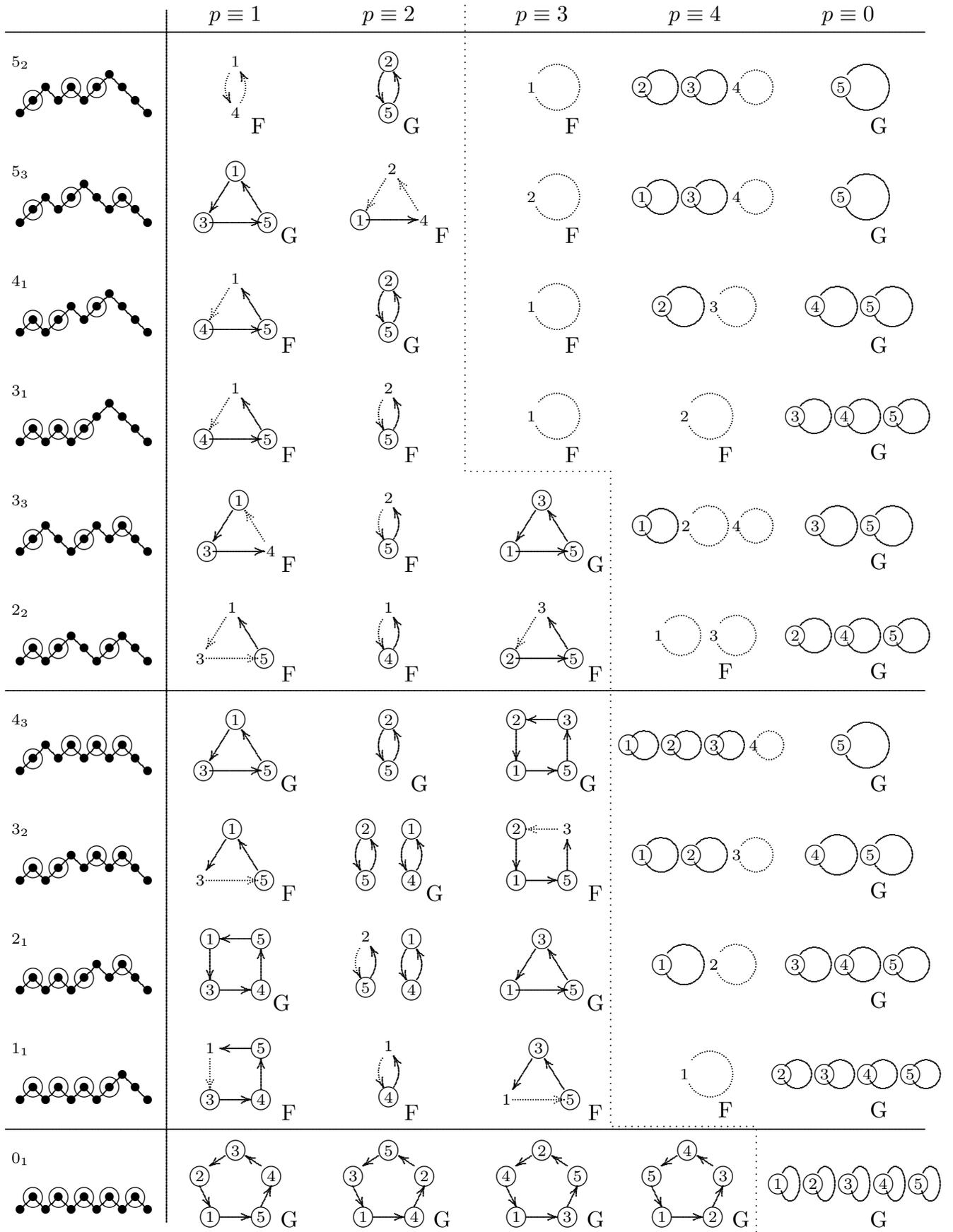
\setdots <1mm>
\plot .58 1  .58 -.95  1.5  -.95  1.5 -9  2.5 -9  2.5 -18.7 /
\setsolid
\put{$p\equiv 1$} at 1 1.3
\put{$p\equiv 2$} at 2 1.3
\put{$p\equiv 3$} at 3 1.3
\put{$p\equiv 4$} at 4 1.3
\put{$p\equiv 0$} at 5 1.3

%==================THE ROOFS===================
\put{\beginpicture
\setcoordinatesystem units <.2cm,.2cm>
\put{$\ssize 10_1$} at 0 4
\multiput{$\bullet$} at 0 0  1 1  2 2  3 3  4 4  5 5  6 4  7 3  8 2  9 1  10 0 /
\plot 0 0  5 5  10 0  /
\put{$\bigcirc$} at 1 1
\endpicture} at 0 0
%========================================
\put{\beginpicture
\setcoordinatesystem units <.2cm,.2cm>
\put{$\ssize 9_1$} at 0 4
\multiput{$\bullet$} at 0 0  1 1  2 2  3 3  4 4  5 3  6 4  7 3  8 2  9 1  10 0 /
\plot 0 0  4 4  5 3  6 4  10 0  /
\multiput{$\bigcirc$} at 1 1  6 4 /
\endpicture} at 0 -2
%========================================
\put{\beginpicture
\setcoordinatesystem units <.2cm,.2cm>
\put{$\ssize 8_1$} at 0 4
\multiput{$\bullet$} at 0 0  1 1  2 2  3 3  4 2  5 3  6 4  7 3  8 2  9 1  10 0 /
\plot 0 0  3 3  4 2  6 4  10 0  /
\multiput{$\bigcirc$} at 1 1  5 3 /
\endpicture} at 0 -4
%========================================
\put{\beginpicture
\setcoordinatesystem units <.2cm,.2cm>
\put{$\ssize 7_1$} at 0 4
\multiput{$\bullet$} at 0 0  1 1  2 2  3 1  4 2  5 3  6 4  7 3  8 2  9 1  10 0 /
\plot 0 0  2 2  3 1  6 4  10 0  /
\multiput{$\bigcirc$} at 1 1  4 2 /
\endpicture} at 0 -6
%========================================
\put{\beginpicture
\setcoordinatesystem units <.2cm,.2cm>
\put{$\ssize 6_1$} at 0 4
\multiput{$\bullet$} at 0 0  1 1  2 0  3 1  4 2  5 3  6 4  7 3  8 2  9 1  10 0 /
\plot 0 0  1 1  2 0  6 4  10 0  /
\multiput{$\bigcirc$} at 1 1  3 1 /
\endpicture} at 0 -8
%========================================
\put{\beginpicture
\setcoordinatesystem units <.2cm,.2cm>
\put{$\ssize 6_3$} at 0 4
\multiput{$\bullet$} at 0 0  1 1  2 2  3 3  4 2  5 1  6 2  7 3  8 2  9 1  10 0 /
\plot 0 0  3 3  5 1  7 3  10 0  /
\multiput{$\bigcirc$} at 1 1  6 2 /
\endpicture} at 0 -10
%========================================
\put{\beginpicture
\setcoordinatesystem units <.2cm,.2cm>
\put{$\ssize 4_2$} at 0 4
\multiput{$\bullet$} at 0 0  1 1  2 2  3 1  4 0  5 1  6 2  7 3  8 2  9 1  10 0 /
\plot 0 0  2 2  4 0  7 3  10 0  /
\multiput{$\bigcirc$} at 1 1  5 1 /
\endpicture} at 0 -12
%========================================
\put{\beginpicture
\setcoordinatesystem units <.2cm,.2cm>
\put{$\ssize 7_2$} at 0 4
\multiput{$\bullet$} at 0 0  1 1  2 2  3 3  4 2  5 3  6 2  7 3  8 2  9 1  10 0 /
\plot 0 0  3 3  4 2  5 3  6 2  7 3  10 0  /
\multiput{$\bigcirc$} at 1 1  5 3  7 3 /
\endpicture} at 0 -14
%========================================
\put{\beginpicture
\setcoordinatesystem units <.2cm,.2cm>
\put{$\ssize 6_2$} at 0 4
\multiput{$\bullet$} at 0 0  1 1  2 2  3 1  4 2  5 3  6 2  7 3  8 2  9 1  10 0 /
\plot 0 0  2 2  3 1  4 2  5 3  6 2  7 3   10 0  /
\multiput{$\bigcirc$} at 1 1  4 2  7 3 /
\endpicture} at 0 -16
%========================================
\put{\beginpicture
\setcoordinatesystem units <.2cm,.2cm>
\put{$\ssize 5_1$} at 0 4
\multiput{$\bullet$} at 0 0  1 1  2 0  3 1  4 2  5 3  6 2  7 3  8 2  9 1  10 0 /
\plot 0 0  1 1  2 0  3 1  5 3  6 2  7 3  10 0  /
\multiput{$\bigcirc$} at 1 1  3 1  7 3 /
\endpicture} at 0 -18
%========================================
%%%%%%%%%%%%%%%%%%%%%%The cycles%%%%%%%%%%%%%%%%%%%%%%%%
%=====================================
\put{\beginpicture
\setcoordinatesystem units <.8cm,.8cm>
\put{$\ssize 5$} at 1 0
\multiput{} at 1.5 0 /
\put{$\bigcirc$} at 1 0
\circulararc 300 degrees from 1.1 -.2  center at 1.5 0 
\endpicture} at 5 0 
%=====================================loop at 1, dotted
\multiput{\beginpicture
\setcoordinatesystem units <.8cm,.8cm>
\put{$\ssize 1$} at 1 0
\multiput{} at 1.5 0 /
\setdots <.4mm>
\circulararc 300 degrees from 1.1 -.2  center at 1.5 0 
\endpicture} at 1 0  2 -4  2 -6  2 -8   3 -12 /
%=====================================loop at 2 dotted
\multiput{\beginpicture
\setcoordinatesystem units <.8cm,.8cm>
\put{$\ssize 2$} at 1 0
\multiput{} at 1.5 0 /
\setdots <.4mm>
\circulararc 300 degrees from 1.1 -.2  center at 1.5 0 
\endpicture} at 2 0  3 -6  3 -8 /
%====================================loop at 3, dotted 
\multiput{\beginpicture
\setcoordinatesystem units <.8cm,.8cm>
\put{$\ssize 3$} at 1 0
\multiput{} at 1.5 0 /
\setdots <.4mm>
\circulararc 300 degrees from 1.1 -.2  center at 1.5 0 
\endpicture} at 3 0  3 -2  4 -8  4  -18 /
%=====================================loop at 4, dotted
\multiput{\beginpicture
\setcoordinatesystem units <.8cm,.8cm>
\put{$\ssize 4$} at 1 0
\multiput{} at 1.5 0 /
\setdots <.4mm>
\circulararc 300 degrees from 1.1 -.2  center at 1.5 0 
\endpicture} at 4 0  4 -2  4 -4 4 -14   / 
%========================================9_1
%=====================================1, a=0
\multiput{\beginpicture
\setcoordinatesystem units <.8cm,.8cm>
\put{$\ssize 5$} at 1 0
\multiput{} at 1.5 0 /
\put{$\bigcirc$} at 1 0
\circulararc 300 degrees from 1.1 -.2  center at 1.5 0 
\endpicture} at 5 -2
                5 -4
                5 -6
                5 -10
                5 -14
                5 -16
                 /
%================================loop at 1, smaller, solid
\multiput{\beginpicture
\setcoordinatesystem units <.7cm,.7cm>
 \put{$\ssize 1$} at 1 0
 \multiput{} at 1.5 0 /
 \put{$\bigcirc$} at 1 0
 \circulararc 300 degrees from 1.1 -.2  center at 1.5 0 
\endpicture} at 1.82 -2  2.82 -16  2.82 -18 /
%================================loop at 2, smaller, solid
\multiput{\beginpicture
\setcoordinatesystem units <.7cm,.7cm>
 \put{$\ssize 2$} at 1 0
 \multiput{} at 1.5 0 /
 \put{$\bigcirc$} at 1 0
 \circulararc 300 degrees from 1.1 -.2  center at 1.5 0 
\endpicture} at 2.82 -4  3.82 -10 /
%================================loop at 3, smaller, solid
\multiput{\beginpicture
\setcoordinatesystem units <.7cm,.7cm>
 \put{$\ssize 3$} at 1 0
 \multiput{} at 1.5 0 /
 \put{$\bigcirc$} at 1 0
 \circulararc 300 degrees from 1.1 -.2  center at 1.5 0 
\endpicture} at 4.82 -12  3.82 -6  3.82 -16 /

%================================loop at 4, smaller, solid
\multiput{\beginpicture
\setcoordinatesystem units <.7cm,.7cm>
 \put{$\ssize 4$} at 1 0
 \multiput{} at 1.5 0 /
 \put{$\bigcirc$} at 1 0
 \circulararc 300 degrees from 1.1 -.2  center at 1.5 0 
\endpicture} at 4.82 -8 4.82 -18 /

%================================loop at 5, smaller, solid
\multiput{\beginpicture
\setcoordinatesystem units <.7cm,.7cm>
 \put{$\ssize 5$} at 1 0
 \multiput{} at 1.5 0 /
 \put{$\bigcirc$} at 1 0
 \circulararc 300 degrees from 1.1 -.2  center at 1.5 0 
\endpicture} at 5.18 -8  5.18 -12 5.18 -18 /

%================================loop at 1, smaller, dotted
\multiput{\beginpicture
\setcoordinatesystem units <.7cm,.7cm>
 \put{$\ssize 1$} at 1 0
 \multiput{} at 1.5 0 /
\setdots <.4mm>
 \circulararc 300 degrees from 1.1 -.2  center at 1.5 0 
\endpicture} at 2.82 -10 /
%================================loop at 2, smaller, dotted
\multiput{\beginpicture
\setcoordinatesystem units <.7cm,.7cm>
 \put{$\ssize 2$} at 1 0
 \multiput{} at 1.5 0 /
\setdots <.4mm>
 \circulararc 300 degrees from 1.1 -.2  center at 1.5 0 
\endpicture} at 2.18 -2  3.82 -12  3.18 -16  3.18 -18 /
%================================loop at 3, smaller, dotted
\multiput{\beginpicture
\setcoordinatesystem units <.7cm,.7cm>
 \put{$\ssize 3$} at 1 0
 \multiput{} at 1.5 0 /
\setdots <.4mm>
 \circulararc 300 degrees from 1.1 -.2  center at 1.5 0 
\endpicture} at 3.18 -4  3.18 -10 /

%================================loop at 4, smaller, dotted
\multiput{\beginpicture
\setcoordinatesystem units <.7cm,.7cm>
 \put{$\ssize 4$} at 1 0
 \multiput{} at 1.5 0 /
\setdots <.4mm>
 \circulararc 300 degrees from 1.1 -.2  center at 1.5 0 
\endpicture} at 4.18 -6  4.18 -10  4.18 -12  4.18 -16 /

%================================loop at 1, small, solid
\multiput{\beginpicture
\setcoordinatesystem units <.6cm,.6cm>
 \put{$\ssize 1$} at 1 0
 \multiput{} at 1.5 0 /
 \put{$\bigcirc$} at 1 0
 \circulararc 300 degrees from 1.1 -.2  center at 1.5 0 
\endpicture} at 2.65 -14 /
%================================loop at 2, small, solid
\multiput{\beginpicture
\setcoordinatesystem units <.6cm,.6cm>
 \put{$\ssize 2$} at 1 0
 \multiput{} at 1.5 0 /
 \put{$\bigcirc$} at 1 0
 \circulararc 300 degrees from 1.1 -.2  center at 1.5 0 
\endpicture} at 3 -14 /
%================================loop at 3, small, dotted
\multiput{\beginpicture
\setcoordinatesystem units <.6cm,.6cm>
 \put{$\ssize 3$} at 1 0
 \multiput{} at 1.5 0 /
\setdots <.4mm>
 \circulararc 300 degrees from 1.1 -.2  center at 1.5 0 
\endpicture} at 3.35 -14 /
%%%%%%%%%%%%%%%%%%%%%%%The 2-cycles%%%%%%%%%%%%%%%%%%%
%================================15solid
\multiput{\beginpicture
\setcoordinatesystem units <.8cm,.8cm>
\put{$\ssize 1$} at 1 0 
\put{$\ssize 5$} at 1 1
\multiput{$\bigcirc$} at 1 0  1 1 /
\arr{0.89 0.22}{0.9 0.2}
\arr{1.11 0.78}{1.1 0.8}
\setquadratic
\plot 0.89 0.2 0.8 0.5 0.9 0.8 /
\plot 1.1 0.2 1.2 0.5 1.1 0.8 /
\setdots <.4mm>
\endpicture} at 1 -2 /

%================================15halfdotted
\multiput{\beginpicture
\setcoordinatesystem units <.8cm,.8cm>
\put{$\ssize 1$} at 1 0 
\put{$\ssize 5$} at 1 1
\multiput{$\bigcirc$} at 1 1  /
\arr{0.89 0.22}{0.9 0.2}
\arr{1.11 0.78}{1.1 0.8}
\setquadratic
\plot 0.89 0.2 0.8 0.5 0.9 0.8 /
\setdots <.4mm>
\plot 1.1 0.2 1.2 0.5 1.1 0.8 /

\endpicture} at 1 -4  1 -6  1 -8 /
%================================25solid
\multiput{\beginpicture
\setcoordinatesystem units <.8cm,.8cm>
\put{$\ssize 5$} at 1 0 
\put{$\ssize 2$} at 1 1
\multiput{$\bigcirc$} at 1 0  1 1 /
\arr{0.89 0.22}{0.9 0.2}
\arr{1.11 0.78}{1.1 0.8}
\setquadratic
\plot 1.1 0.2 1.2 0.5 1.1 0.8 /
\plot 0.89 0.2 0.8 0.5 0.9 0.8 /
\endpicture} at 2 -10 /
%================================25halfdotted
\multiput{\beginpicture
\setcoordinatesystem units <.8cm,.8cm>
\put{$\ssize 5$} at 1 0 
\put{$\ssize 2$} at 1 1
\multiput{$\bigcirc$} at 1 0   /
\arr{0.89 0.22}{0.9 0.2}
\arr{1.11 0.78}{1.1 0.8}
\setquadratic
\plot 1.1 0.2 1.2 0.5 1.1 0.8 /
\setdots <.4mm>
\plot 0.89 0.2 0.8 0.5 0.9 0.8 /
\endpicture} at 2 -12 /

%================================14halfdotted
\multiput{\beginpicture
\setcoordinatesystem units <.8cm,.8cm>
\put{$\ssize 4$} at 1 0 
\put{$\ssize 1$} at 1 1
\multiput{$\bigcirc$} at  1 1 /
\arr{0.89 0.22}{0.9 0.2}
\arr{1.11 0.78}{1.1 0.8}
\setquadratic
\plot 0.89 0.2 0.8 0.5 0.9 0.8 /
\setdots <.4mm>
\plot 1.1 0.2 1.2 0.5 1.1 0.8 /

\endpicture} at 1 -14  1 -16 /
%================================14dotted
\multiput{\beginpicture
\setcoordinatesystem units <.8cm,.8cm>
\put{$\ssize 4$} at 1 0 
\put{$\ssize 1$} at 1 1
\arr{0.89 0.22}{0.9 0.2}
\arr{1.11 0.78}{1.1 0.8}
\setquadratic
\setdots <.4mm>
\plot 1.1 0.2 1.2 0.5 1.1 0.8 /
\plot 0.89 0.2 0.8 0.5 0.9 0.8 /

\endpicture} at 1 -10  1 -12 /
%%%%%%%%%%%%%%%%%3-Zyklen%%%%%%%%%%%%%%
%================================125solid
\multiput{\beginpicture
\setcoordinatesystem units <.5cm,.8cm>
\put{$\ssize 1$} at 0 0 
\put{$\ssize 2$} at 1 1
\put{$\ssize 5$} at 2 0
\multiput{$\bigcirc$} at 0 0  1 1  2 0  /
\arr{0.2 0}{1.8 0}
\arr{1.8 0.2}{1.2 0.8}
\arr{0.8 0.8}{0.2 0.2}
\endpicture} at 2 -14 /
%================================145solid
\multiput{\beginpicture
\setcoordinatesystem units <.5cm,.8cm>
\put{$\ssize 4$} at 0 0 
\put{$\ssize 1$} at 1 1
\put{$\ssize 5$} at 2 0
\multiput{$\bigcirc$} at 0 0  1 1  2 0  /
\arr{0.2 0}{1.8 0}
\arr{1.8 0.2}{1.2 0.8}
\arr{0.8 0.8}{0.2 0.2}
\endpicture} at 1 -18 /
%================================125third
\multiput{\beginpicture
\setcoordinatesystem units <.5cm,.8cm>
\put{$\ssize 1$} at 0 0 
\put{$\ssize 2$} at 1 1
\put{$\ssize 5$} at 2 0
\multiput{$\bigcirc$} at 0 0   2 0  /
\arr{0.2 0}{1.8 0}
\arr{1.8 0.2}{1.2 0.8}
\setdots <.4mm>
\arr{0.8 0.8}{0.2 0.2}
\endpicture} at 2 -16  2 -18 /

\multiput{G} at 
 1.15 -2.7
 2.15 -10.7
 2.35 -14.7
 1.35 -18.7
 5.2 -0.7 
 5.2 -2.7   
 5.2 -4.7     
 5.2 -6.7  
 5.2 -8.65  
 5.2 -10.7   
 5.2 -12.65   
 5.2 -14.7  
 5.2 -16.7  
 5.2 -18.65   /
\multiput{F} at
 1.2 -0.7 
 2.2 -0.7
 3.2 -0.7
 4.2 -0.7
 3.2 -2.7
 4.2 -2.7
 1.15 -4.7 
 2.2 -4.7
 4.2 -4.7
 1.15 -6.7 
 2.2 -6.7
 3.2 -6.7
 1.15 -8.7 
 2.2 -8.7
 3.2 -8.7
 4.2 -8.7
 1.15 -10.7
 3.2 -10.65
 1.15 -12.7 
 2.15 -12.7 
 3.2 -12.7
 4.2 -12.65
 1.15 -14.7
 4.2 -14.7
 1.15 -16.7
 2.35 -16.7
 2.35 -18.7
 4.2 -18.7
/

\endpicture}
$$

$$
{\beginpicture
%\setcoordinatesystem units <2.4cm,.94cm>
\setcoordinatesystem units <2.4cm,.86cm>
%========================================second page
\plot -.5 1  5.5 1 /
\plot -.5 -11  5.5 -11 /
\plot -.5 -19  5.5 -19 /
\plot 0.55 1.4  0.55 -20.7 /

\setdots <1mm>
\plot 2.5 1.5  2.5 -7  3.45 -7  3.45 -18.95  4.4 -18.95  4.4 -21  /
\setsolid

\put{$p\equiv 1$} at 1 1.3
\put{$p\equiv 2$} at 2 1.3
\put{$p\equiv 3$} at 3 1.3
\put{$p\equiv 4$} at 4 1.3
\put{$p\equiv 0$} at 5 1.3

%==================THE ROOFS===================
\put{\beginpicture
\setcoordinatesystem units <.2cm,.2cm>
\put{$\ssize 5_2$} at 0 4
\multiput{$\bullet$} at 0 0  1 1  2 2  3 1  4 2  5 1  6 2  7 3  8 2  9 1  10 0 /
\plot 0 0  2 2  3 1  4 2  5 1  6 2  7 3  10 0  /
\multiput{$\bigcirc$} at 1 1  4 2  6 2  /
\endpicture} at 0 0
%========================================
\put{\beginpicture
\setcoordinatesystem units <.2cm,.2cm>
\put{$\ssize 5_3$} at 0 4
\multiput{$\bullet$} at 0 0  1 1  2 2  3 1  4 2  5 3  6 2  7 1  8 2  9 1  10 0 /
\plot 0 0  1 1  2 2  3 1  4 2  5 3  6 2  7 1  8 2  9 1  10 0  /
\multiput{$\bigcirc$} at 1 1  4 2  8 2   /
\endpicture} at 0 -2
%========================================
\put{\beginpicture
\setcoordinatesystem units <.2cm,.2cm>
\put{$\ssize 4_1$} at 0 4
\multiput{$\bullet$} at 0 0  1 1  2 0  3 1  4 2  5 1  6 2  7 3  8 2  9 1  10 0 /
\plot 0 0  1 1  2 0  3 1  4 2  5 1  6 2  7 3  8 2  9 1  10 0  /
\multiput{$\bigcirc$} at 1 1  3 1  6 2 /
\endpicture} at 0 -4
%========================================

\put{\beginpicture
\setcoordinatesystem units <.2cm,.2cm>
\put{$\ssize 3_1$} at 0 4
\multiput{$\bullet$} at 0 0  1 1  2 0  3 1  4 0  5 1  6 2  7 3  8 2  9 1  10 0 /
\plot 0 0  1 1  2 0  3 1  4 0  5 1  6 2  7 3  8 2  9 1  10 0  /
\multiput{$\bigcirc$} at 1 1  3 1  5 1 /
\endpicture} at 0 -6
%========================================
\put{\beginpicture
\setcoordinatesystem units <.2cm,.2cm>
\put{$\ssize 3_3$} at 0 4
\multiput{$\bullet$} at 0 0  1 1  2 2  3 1  4 0  5 1  6 2  7 1  8 2  9 1  10 0 /
\plot 0 0  1 1  2 2  3 1  4 0  5 1  6 2  7 1  8 2  9 1  10 0  /
\multiput{$\bigcirc$} at 1 1  5 1  8 2 /
\endpicture} at 0 -8
%========================================
\put{\beginpicture
\setcoordinatesystem units <.2cm,.2cm>
\put{$\ssize 2_2$} at 0 4
\multiput{$\bullet$} at 0 0  1 1  2 0  3 1  4 2  5 1  6 0  7 1  8 2  9 1  10 0 /
\plot 0 0  1 1  2 0  3 1  4 2  5 1  6 0  7 1  8 2  9 1  10 0 /
\multiput{$\bigcirc$} at 1 1  3 1  7 1 /
\endpicture} at 0 -10
%========================================
%========================================
\put{\beginpicture
\setcoordinatesystem units <.2cm,.2cm>
\put{$\ssize 4_3$} at 0 4
\multiput{$\bullet$} at 0 0  1 1  2 2  3 1  4 2  5 1  6 2  7 1  8 2  9 1  10 0 /
\plot 0 0  1 1  2 2  3 1  4 2  5 1  6 2  7 1  8 2  9 1  10 0  /
\multiput{$\bigcirc$} at 1 1  4 2  6 2  8 2  /
\endpicture} at 0 -12
%========================================
\put{\beginpicture
\setcoordinatesystem units <.2cm,.2cm>
\put{$\ssize 3_2$} at 0 4
\multiput{$\bullet$} at 0 0  1 1  2 0  3 1  4 2  5 1  6 2  7 1  8 2  9 1  10 0 /
\plot 0 0  1 1  2 0  3 1  4 2  5 1  6 2  7 1  8 2  9 1  10 0 /
\multiput{$\bigcirc$} at 1 1  3 1  6 2  8 2 /
\endpicture} at 0 -14
%========================================
\put{\beginpicture
\setcoordinatesystem units <.2cm,.2cm>
\put{$\ssize 2_1$} at 0 4
\multiput{$\bullet$} at 0 0  1 1  2 0  3 1  4 0  5 1  6 2  7 1  8 2  9 1  10 0 /
\plot 0 0  1 1  2 0  3 1  4 0  5 1  6 2  7 1  8 2  9 1  10 0 /
\multiput{$\bigcirc$} at 1 1  3 1  5 1  8 2 /
\endpicture} at 0 -16
%========================================
\put{\beginpicture
\setcoordinatesystem units <.2cm,.2cm>
\put{$\ssize 1_1$} at 0 4
\multiput{$\bullet$} at 0 0  1 1  2 0  3 1  4 0  5 1  6 0  7 1  8 2  9 1  10 0 /
\plot 0 0  1 1  2 0  3 1  4 0  5 1  6 0  7 1  8 2  9 1  10 0 /
\multiput{$\bigcirc$} at 1 1  3 1  5 1  7 1  /
\endpicture} at 0 -18
%========================================
%========================================
\put{\beginpicture
\setcoordinatesystem units <.2cm,.2cm>
\put{$\ssize 0_1$} at 0 4
\multiput{$\bullet$} at 0 0  1 1  2 0  3 1  4 0  5 1  6 0  7 1   9 1  10 0  8 0 /
\plot 0 0  1 1  2 0  3 1  4 0  5 1  6 0  7 1  8 0  9 1  10 0 /
\multiput{$\bigcirc$} at 1 1  3 1  5 1  7 1  9 1  /
\endpicture} at 0 -20
%========================================

%%%%%%%%%%%%%%%%%%%%%%The cycles%%%%%%%%%%%%%%%%%%%%%%%%
%=====================================loop at 1, dotted
\multiput{\beginpicture
\setcoordinatesystem units <.8cm,.8cm>
\put{$\ssize 1$} at 1 0
\multiput{} at 1.5 0 /
\setdots <.4mm>
\circulararc 300 degrees from 1.1 -.2  center at 1.5 0 
\endpicture} at 3 0   3 -4  3 -6    4 -18 /
%=====================================loop at 2 dotted
\multiput{\beginpicture
\setcoordinatesystem units <.8cm,.8cm>
\put{$\ssize 2$} at 1 0
\multiput{} at 1.5 0 /
\setdots <.4mm>
\circulararc 300 degrees from 1.1 -.2  center at 1.5 0 
\endpicture} at 3 -2   4 -6   /

%=====================================loop at 5
\multiput{\beginpicture
\setcoordinatesystem units <.8cm,.8cm>
\put{$\ssize 5$} at 1 0
\multiput{} at 1.5 0 /
\put{$\bigcirc$} at 1 0
\circulararc 300 degrees from 1.1 -.2  center at 1.5 0 
\endpicture} at 5 -0
                5 -2
                5 -12   /

%================================loop at 1, smaller, solid
\multiput{\beginpicture
\setcoordinatesystem units <.7cm,.7cm>
 \put{$\ssize 1$} at 1 0
 \multiput{} at 1.5 0 /
 \put{$\bigcirc$} at 1 0
 \circulararc 300 degrees from 1.1 -.2  center at 1.5 0 
\endpicture} at 3.82 -16 /
%================================loop at 2, smaller, solid
\multiput{\beginpicture
\setcoordinatesystem units <.7cm,.7cm>
 \put{$\ssize 2$} at 1 0
 \multiput{} at 1.5 0 /
 \put{$\bigcirc$} at 1 0
 \circulararc 300 degrees from 1.1 -.2  center at 1.5 0 
\endpicture} at 3.82 -4 /
%================================loop at 3, smaller, solid
\multiput{\beginpicture
\setcoordinatesystem units <.7cm,.7cm>
 \put{$\ssize 3$} at 1 0
 \multiput{} at 1.5 0 /
 \put{$\bigcirc$} at 1 0
 \circulararc 300 degrees from 1.1 -.2  center at 1.5 0 
\endpicture} at  4.82  -8 /

%================================loop at 4, smaller, solid
\multiput{\beginpicture
\setcoordinatesystem units <.7cm,.7cm>
 \put{$\ssize 4$} at 1 0
 \multiput{} at 1.5 0 /
 \put{$\bigcirc$} at 1 0
 \circulararc 300 degrees from 1.1 -.2  center at 1.5 0 
\endpicture} at 4.82 -4  4.82 -14 /

%================================loop at 5, smaller, solid
\multiput{\beginpicture
\setcoordinatesystem units <.7cm,.7cm>
 \put{$\ssize 5$} at 1 0
 \multiput{} at 1.5 0 /
 \put{$\bigcirc$} at 1 0
 \circulararc 300 degrees from 1.1 -.2  center at 1.5 0 
\endpicture} at 5.18 -4  5.18 -8  5.18 -14 /

%================================loop at 1, smaller, dotted
\multiput{\beginpicture
\setcoordinatesystem units <.7cm,.7cm>
 \put{$\ssize 1$} at 1 0
 \multiput{} at 1.5 0 /
\setdots <.4mm>
 \circulararc 300 degrees from 1.1 -.2  center at 1.5 0 
\endpicture} at 3.82 -10 /
%================================loop at 2, smaller, dotted
\multiput{\beginpicture
\setcoordinatesystem units <.7cm,.7cm>
 \put{$\ssize 2$} at 1 0
 \multiput{} at 1.5 0 /
\setdots <.4mm>
 \circulararc 300 degrees from 1.1 -.2  center at 1.5 0 
\endpicture} at 4.18 -16  4 -8 /
%================================loop at 3, smaller, dotted
\multiput{\beginpicture
\setcoordinatesystem units <.7cm,.7cm>
 \put{$\ssize 3$} at 1 0
 \multiput{} at 1.5 0 /
\setdots <.4mm>
 \circulararc 300 degrees from 1.1 -.2  center at 1.5 0 
\endpicture} at 4.18 -4  4.18 -10 /

%================================loop at 1, small, sold
\multiput{\beginpicture
\setcoordinatesystem units <.6cm,.6cm>
 \put{$\ssize 1$} at 1 0
 \multiput{} at 1.5 0 /
 \put{$\bigcirc$} at 1 0
 \circulararc 300 degrees from 1.1 -.2  center at 1.5 0 
\endpicture} at 3.68 -2  3.68 -8  3.68 -14 /
%================================loop at 2, small, solid
\multiput{\beginpicture
\setcoordinatesystem units <.6cm,.6cm>
 \put{$\ssize 2$} at 1 0
 \multiput{} at 1.5 0 /
 \put{$\bigcirc$} at 1 0
 \circulararc 300 degrees from 1.1 -.2  center at 1.5 0 
\endpicture} at 3.68 0  4 -14  4.68 -10 /
%================================loop at 3, small, solid
\multiput{\beginpicture
\setcoordinatesystem units <.6cm,.6cm>
 \put{$\ssize 3$} at 1 0
 \multiput{} at 1.5 0 /
 \put{$\bigcirc$} at 1 0
 \circulararc 300 degrees from 1.1 -.2  center at 1.5 0 
\endpicture} at   4.68 -6  4.68 -16   4 0  4 -2 /
%================================loop at 4, small, solid
\multiput{\beginpicture
\setcoordinatesystem units <.6cm,.6cm>
 \put{$\ssize 4$} at 1 0
 \multiput{} at 1.5 0 /
 \put{$\bigcirc$} at 1 0
 \circulararc 300 degrees from 1.1 -.2  center at 1.5 0 
\endpicture} at   5 -6  5 -10  5 -16 /
%================================loop at 5, small, solid
\multiput{\beginpicture
\setcoordinatesystem units <.6cm,.6cm>
 \put{$\ssize 5$} at 1 0
 \multiput{} at 1.5 0 /
 \put{$\bigcirc$} at 1 0
 \circulararc 300 degrees from 1.1 -.2  center at 1.5 0 
\endpicture} at   5.32 -6  5.32 -10  5.32 -16 /

%================================loop at 3, small, dotted
\multiput{\beginpicture
\setcoordinatesystem units <.6cm,.6cm>
 \put{$\ssize 3$} at 1 0
 \multiput{} at 1.5 0 /
\setdots <.4mm>
 \circulararc 300 degrees from 1.1 -.2  center at 1.5 0 
\endpicture} at   4.32 -14 /
%================================loop at 4, small, dotted
\multiput{\beginpicture
\setcoordinatesystem units <.6cm,.6cm>
 \put{$\ssize 4$} at 1 0
 \multiput{} at 1.5 0 /
\setdots <.4mm>
 \circulararc 300 degrees from 1.1 -.2  center at 1.5 0 
\endpicture} at 4.32 0  4.32 -2  4.32 -8 /

%================================loop at 1, very small, solid
\multiput{\beginpicture
\setcoordinatesystem units <.5cm,.5cm>
 \put{$\ssize 1$} at 1 0
 \multiput{} at 1.5 0 /
 \put{$\bigcirc$} at 1 0
 \circulararc 300 degrees from 1.1 -.2  center at 1.5 0 
\endpicture} at  3.58 -12   /
%================================loop at 2, very small, solid
\multiput{\beginpicture
\setcoordinatesystem units <.5cm,.5cm>
 \put{$\ssize 2$} at 1 0
 \multiput{} at 1.5 0 /
 \put{$\bigcirc$} at 1 0
 \circulararc 300 degrees from 1.1 -.2  center at 1.5 0 
\endpicture} at   3.86 -12  4.58 -18 /
%================================loop at 3, very small, solid
\multiput{\beginpicture
\setcoordinatesystem units <.5cm,.5cm>
 \put{$\ssize 3$} at 1 0
 \multiput{} at 1.5 0 /
 \put{$\bigcirc$} at 1 0
 \circulararc 300 degrees from 1.1 -.2  center at 1.5 0 
\endpicture} at   4.86 -18  4.14 -12 /
%================================loop at 4, very small, solid
\multiput{\beginpicture
\setcoordinatesystem units <.5cm,.5cm>
 \put{$\ssize 4$} at 1 0
 \multiput{} at 1.5 0 /
 \put{$\bigcirc$} at 1 0
 \circulararc 300 degrees from 1.1 -.2  center at 1.5 0 
\endpicture} at   5.14 -18  /
%================================loop at 5, very small, solid
\multiput{\beginpicture
\setcoordinatesystem units <.5cm,.5cm>
 \put{$\ssize 5$} at 1 0
 \multiput{} at 1.5 0 /
 \put{$\bigcirc$} at 1 0
 \circulararc 300 degrees from 1.1 -.2  center at 1.5 0 
\endpicture} at  5.42 -18 /
%================================loop at 4, very small, dotted
\multiput{\beginpicture
\setcoordinatesystem units <.5cm,.5cm>
 \put{$\ssize 4$} at 1 0
 \multiput{} at 1.5 0 /
\setdots <.4mm>
 \circulararc 300 degrees from 1.1 -.2  center at 1.5 0 
\endpicture} at  4.42 -12 /

%%%%%%%%%%%%%%%%%%52cycles%%%%%%%%%%%%%%%%%%%%%%
%================================25solid
\multiput{\beginpicture
\setcoordinatesystem units <.8cm,.8cm>
\put{$\ssize 5$} at 1 0 
\put{$\ssize 2$} at 1 1
\multiput{$\bigcirc$} at 1 0  1 1 /
\arr{0.89 0.22}{0.9 0.2}
\arr{1.11 0.78}{1.1 0.8}
\setquadratic
\plot 1.1 0.2 1.2 0.5 1.1 0.8 /
\plot 0.89 0.2 0.8 0.5 0.9 0.8 /
\endpicture} at 2 0  2 -4  2 -12  1.85 -14  /
%================================25halfdotted
\multiput{\beginpicture
\setcoordinatesystem units <.8cm,.8cm>
\put{$\ssize 5$} at 1 0 
\put{$\ssize 2$} at 1 1
\multiput{$\bigcirc$} at 1 0   /
\arr{0.89 0.22}{0.9 0.2}
\arr{1.11 0.78}{1.1 0.8}
\setquadratic
\plot 1.1 0.2 1.2 0.5 1.1 0.8 /
\setdots <.4mm>
\plot 0.89 0.2 0.8 0.5 0.9 0.8 /
\endpicture} at 2 -6  2 -8  1.85 -16  /

%================================14halfdotted
\multiput{\beginpicture
\setcoordinatesystem units <.8cm,.8cm>
\put{$\ssize 4$} at 1 0 
\put{$\ssize 1$} at 1 1
\multiput{$\bigcirc$} at  1 0 /
\arr{0.89 0.22}{0.9 0.2}
\arr{1.11 0.78}{1.1 0.8}
\setquadratic
\plot 1.1 0.2 1.2 0.5 1.1 0.8 /
\setdots <.4mm>
\plot 0.89 0.2 0.8 0.5 0.9 0.8 /

\endpicture} at 2 -10  2 -18   /
%================================14dotted
\multiput{\beginpicture
\setcoordinatesystem units <.8cm,.8cm>
\put{$\ssize 4$} at 1 0 
\put{$\ssize 1$} at 1 1
\arr{0.89 0.22}{0.9 0.2}
\arr{1.11 0.78}{1.1 0.8}
\setquadratic
\setdots <.4mm>
\plot 1.1 0.2 1.2 0.5 1.1 0.8 /
\plot 0.89 0.2 0.8 0.5 0.9 0.8 /

\endpicture} at 1 -0  /
%================================14solid
\multiput{\beginpicture
\setcoordinatesystem units <.8cm,.8cm>
\put{$\ssize 4$} at 1 0 
\put{$\ssize 1$} at 1 1
\multiput{$\bigcirc$} at 1 0  1 1 /
\arr{0.89 0.22}{0.9 0.2}
\arr{1.11 0.78}{1.1 0.8}
\setquadratic
\plot 1.1 0.2 1.2 0.5 1.1 0.8 /
\plot 0.89 0.2 0.8 0.5 0.9 0.8 /

\endpicture} at 2.15 -14  2.15 -16 /

%%%%%%%%%%%%%%%%%3-Zyklen%%%%%%%%%%%%%%
%================================135solid
\multiput{\beginpicture
\setcoordinatesystem units <.5cm,.8cm>
\put{$\ssize 3$} at 0 0 
\put{$\ssize 1$} at 1 1
\put{$\ssize 5$} at 2 0
\multiput{$\bigcirc$} at 0 0  1 1  2 0  /
\arr{0.2 0}{1.8 0}
\arr{1.8 0.2}{1.2 0.8}
\arr{0.8 0.8}{0.2 0.2}
\endpicture} at 1 -2  1 -12 /
%================================135third
\multiput{\beginpicture
\setcoordinatesystem units <.5cm,.8cm>
\put{$\ssize 3$} at 0 0 
\put{$\ssize 1$} at 1 1
\put{$\ssize 5$} at 2 0
\multiput{$\bigcirc$} at  1 1  2 0  /
\arr{1.8 0.2}{1.2 0.8}
\arr{0.8 0.8}{0.2 0.2}
\setdots <.4mm>
\arr{0.2 0}{1.8 0}
\endpicture} at 1 -14 /
%================================145solid
\multiput{\beginpicture
\setcoordinatesystem units <.5cm,.8cm>
\put{$\ssize 4$} at 0 0 
\put{$\ssize 1$} at 1 1
\put{$\ssize 5$} at 2 0
\multiput{$\bigcirc$} at 0 0  1 1  2 0  /
\arr{0.2 0}{1.8 0}
\arr{1.8 0.2}{1.2 0.8}
\arr{0.8 0.8}{0.2 0.2}
\endpicture} at  /
%================================145third 
\multiput{\beginpicture
\setcoordinatesystem units <.5cm,.8cm>
\put{$\ssize 4$} at 0 0 
\put{$\ssize 1$} at 1 1
\put{$\ssize 5$} at 2 0
\multiput{$\bigcirc$} at  0 0   2 0  /
\arr{0.2 0}{1.8 0}
\arr{1.8 0.2}{1.2 0.8}
\setdots <.4mm>
\arr{0.8 0.8}{0.2 0.2}
\endpicture} at 1 -4  1 -6 /
%================================134third
\multiput{\beginpicture
\setcoordinatesystem units <.5cm,.8cm>
\put{$\ssize 3$} at 0 0 
\put{$\ssize 1$} at 1 1
\put{$\ssize 4$} at 2 0
\multiput{$\bigcirc$} at 0 0  1 1   /
\arr{0.2 0}{1.8 0}
\arr{0.8 0.8}{0.2 0.2}
\setdots <.4mm>
\arr{1.8 0.2}{1.2 0.8}
\endpicture} at 1 -8 /
%================================135twothird
\multiput{\beginpicture
\setcoordinatesystem units <.5cm,.8cm>
\put{$\ssize 3$} at 0 0 
\put{$\ssize 1$} at 1 1
\put{$\ssize 5$} at 2 0
\multiput{$\bigcirc$} at  2 0  /
\arr{1.8 0.2}{1.2 0.8}
\setdots <.4mm>
\arr{0.2 0}{1.8 0}
\arr{0.8 0.8}{0.2 0.2}
\endpicture} at 1 -10 /

%================================142twothird
\multiput{\beginpicture
\setcoordinatesystem units <.5cm,.8cm>
\put{$\ssize 1$} at 0 0 
\put{$\ssize 2$} at 1 1
\put{$\ssize 4$} at 2 0
\multiput{$\bigcirc$} at  0 0  /
\arr{0.2 0}{1.8 0}
\setdots <.4mm>
\arr{1.8 0.2}{1.2 0.8}
\arr{0.8 0.8}{0.2 0.2}
\endpicture} at 2 -2 /

%================================153
\multiput{\beginpicture
\setcoordinatesystem units <.5cm,.8cm>
\put{$\ssize 1$} at 0 0 
\put{$\ssize 3$} at 1 1
\put{$\ssize 5$} at 2 0
\multiput{$\bigcirc$} at 0 0  1 1  2 0  /
\arr{1.8 0.2}{1.2 0.8}
\arr{0.2 0}{1.8 0}
\arr{0.8 0.8}{0.2 0.2}
\setdots <.4mm>
\endpicture} at 3 -8  3 -16 /
%================================153third
\multiput{\beginpicture
\setcoordinatesystem units <.5cm,.8cm>
\put{$\ssize 1$} at 0 0 
\put{$\ssize 3$} at 1 1
\put{$\ssize 5$} at 2 0
\multiput{$\bigcirc$} at  1 1  2 0  /
\arr{1.8 0.2}{1.2 0.8}
\arr{0.8 0.8}{0.2 0.2}
\setdots <.4mm>
\arr{0.2 0}{1.8 0}
\endpicture} at 3 -18 /
%================================253third
\multiput{\beginpicture
\setcoordinatesystem units <.5cm,.8cm>
\put{$\ssize 2$} at 0 0 
\put{$\ssize 3$} at 1 1
\put{$\ssize 5$} at 2 0
\multiput{$\bigcirc$} at 0 0   2 0  /
\arr{1.8 0.2}{1.2 0.8}
\arr{0.2 0}{1.8 0}
\setdots <.4mm>
\arr{0.8 0.8}{0.2 0.2}
\endpicture} at 3 -10 /

%%%%%%%%%%%%%%%Vierer%%%%%%%%%%%%%%%%%%%
%================================1345
\multiput{\beginpicture
\setcoordinatesystem units <.8cm,.8cm>
\put{$\ssize 1$} at 0 1
\put{$\ssize 3$} at 0 0 
\put{$\ssize 4$} at 1 0
\put{$\ssize 5$} at 1 1
\multiput{$\bigcirc$} at 0 0  1 0  1 1  0 1 /
\arr{0 0.8}{0 0.2}
\arr{0.2 0}{0.8 0}
\arr{0.8 1}{0.2 1}
\arr{1 0.2}{1 0.8}
\endpicture} at 1 -16 /
%================================1345fourth
\multiput{\beginpicture
\setcoordinatesystem units <.8cm,.8cm>
\put{$\ssize 1$} at 0 1
\put{$\ssize 3$} at 0 0 
\put{$\ssize 4$} at 1 0
\put{$\ssize 5$} at 1 1
\multiput{$\bigcirc$} at 0 0  1 0  1 1  /
\setdots <.4mm>
\arr{0 0.8}{0 0.2}
\setsolid
\arr{0.2 0}{0.8 0}
\arr{0.8 1}{0.2 1}
\arr{1 0.2}{1 0.8}
\endpicture} at 1 -18 /
%================================1532
\multiput{\beginpicture
\setcoordinatesystem units <.8cm,.8cm>
\put{$\ssize 2$} at 0 1
\put{$\ssize 1$} at 0 0 
\put{$\ssize 5$} at 1 0
\put{$\ssize 3$} at 1 1
\multiput{$\bigcirc$} at 0 0  1 0  1 1  0 1 /
\arr{0 0.8}{0 0.2}
\arr{0.2 0}{0.8 0}
\arr{0.8 1}{0.2 1}
\arr{1 0.2}{1 0.8}
\endpicture} at 3 -12 /
%================================1235fourth
\multiput{\beginpicture
\setcoordinatesystem units <.8cm,.8cm>
\put{$\ssize 2$} at 0 1
\put{$\ssize 1$} at 0 0 
\put{$\ssize 5$} at 1 0
\put{$\ssize 3$} at 1 1
\multiput{$\bigcirc$} at 0 0  0 1  1 0  /
\arr{0 0.8}{0 0.2}
\arr{0.2 0}{0.8 0}
\setdots <.4mm>
\arr{0.8 1}{0.2 1}
\setsolid
\arr{1 0.2}{1 0.8}
\endpicture} at 3 -14 /

\multiput{G} at 
 2.15 -0.7
 1.35 -2.7
 2.15 -4.7 
 3.35 -8.7
 1.35 -12.7
 2.2 -12.7
 3.3 -12.7
 2.3 -14.7
 1.3 -16.7
 3.35 -16.7
 5.2 -0.7 
 5.2 -2.7   
 5.2 -4.7   
 5.2 -6.65   
 5.2 -8.65  
 5.2 -10.65   
 5.2 -12.7   
 5.2 -14.65   
 5.2 -16.65  
 5.2 -18.65  
 1.35 -20.65
 2.35 -20.65
 3.35 -20.65
 4.3 -20.65
 5.2 -20.65
 /
\multiput{F} at
 1.15 -0.7 
 3.2 -0.7
 2.35 -2.7
 3.2 -2.7
 1.35 -4.7
 3.2 -4.7
 1.35 -6.7
 2.15 -6.7
 3.2 -6.7
 4.2 -6.7
 1.35 -8.7
 2.15 -8.7
 1.35 -10.7
 2.15 -10.7
 3.35 -10.7
 4.2 -10.65
 1.35 -14.7
 3.35 -14.7
1.35 -18.7
 2.15 -18.7
 3.35 -18.7
 4.2 -18.65 
/

%================================loops solid
\multiput{\beginpicture
\setcoordinatesystem units <.55cm,.55cm>
 \put{$\ssize 1$} at 1 0
 \put{$\ssize 2$} at 2 0
 \put{$\ssize 3$} at 3 0
 \put{$\ssize 4$} at 4 0
 \put{$\ssize 5$} at 5 0
 \multiput{$\bigcirc$} at 1 0  2 0  3 0  4 0  5 0 /
 \ellipticalarc axes ratio 1:1.8 300 degrees from 1.1 -.2  center at 1.35 0 
 \ellipticalarc axes ratio 1:1.8 300 degrees from 2.1 -.2  center at 2.35 0 
 \ellipticalarc axes ratio 1:1.8 300 degrees from 3.1 -.2  center at 3.35 0 
 \ellipticalarc axes ratio 1:1.8 300 degrees from 4.1 -.2  center at 4.35 0 
 \ellipticalarc axes ratio 1:1.8 300 degrees from 5.1 -.2  center at 5.35 0 
\endpicture} at 5 -20 /

%================================12345
\put{\beginpicture
\setcoordinatesystem units <.8cm,.8cm>
\put{$\ssize 1$} at 0 0 
\put{$\ssize 2$} at -.2 .8
\put{$\ssize 3$} at 0.5 1.3
\put{$\ssize 4$} at 1.2 .8
\put{$\ssize 5$} at 1 0
\multiput{$\bigcirc$} at 0 0  1 0  -.2 .8  0.5 1.3  1.2 .8 /
\arr{0.2 0}{0.8 0}
\arr{-.2 .6}{0 .2}
\arr{1 .2}{1.2 .6}
\arr{0.3 1.2}{0 1}
\arr{1 1}{0.7 1.2}
\endpicture} at 1 -20

%================================13...
\put{\beginpicture
\setcoordinatesystem units <.8cm,.8cm>
\put{$\ssize 1$} at 0 0 
\put{$\ssize 3$} at -.2 .8
\put{$\ssize 5$} at 0.5 1.3
\put{$\ssize 2$} at 1.2 .8
\put{$\ssize 4$} at 1 0
\multiput{$\bigcirc$} at 0 0  1 0  -.2 .8  0.5 1.3  1.2 .8 /
\arr{0.2 0}{0.8 0}
\arr{-.2 .6}{0 .2}
\arr{1 .2}{1.2 .6}
\arr{0.3 1.2}{0 1}
\arr{1 1}{0.7 1.2}
\endpicture} at 2 -20

%================================12345
\put{\beginpicture
\setcoordinatesystem units <.8cm,.8cm>
\put{$\ssize 1$} at 0 0 
\put{$\ssize 4$} at -.2 .8
\put{$\ssize 2$} at 0.5 1.3
\put{$\ssize 5$} at 1.2 .8
\put{$\ssize 3$} at 1 0
\multiput{$\bigcirc$} at 0 0  1 0  -.2 .8  0.5 1.3  1.2 .8 /
\arr{0.2 0}{0.8 0}
\arr{-.2 .6}{0 .2}
\arr{1 .2}{1.2 .6}
\arr{0.3 1.2}{0 1}
\arr{1 1}{0.7 1.2}
\endpicture} at 3 -20

%================================12345
\put{\beginpicture
\setcoordinatesystem units <.8cm,.8cm>
\put{$\ssize 1$} at 0 0 
\put{$\ssize 5$} at -.2 .8
\put{$\ssize 4$} at 0.5 1.3
\put{$\ssize 3$} at 1.2 .8
\put{$\ssize 2$} at 1 0
\multiput{$\bigcirc$} at 0 0  1 0  -.2 .8  0.5 1.3  1.2 .8 /
\arr{0.2 0}{0.8 0}
\arr{-.2 .6}{0 .2}
\arr{1 .2}{1.2 .6}
\arr{0.3 1.2}{0 1}
\arr{1 1}{0.7 1.2}
\endpicture} at 3.95 -20

\endpicture}
$$
	\bigskip\bigskip
%=======================================================================
{\bf 7.~Nakayama algebras with $p \le s.$}
     \medskip
Some of our results need the condition $p > s$. Here are examples
concerning the possible behavior of algebras with $p\le s$ (see also [CY]).
	   \medskip
{\bf Example 1.} A Nakayama algebra with 
a black cycle (and no other cycles) in the resolution quiver which is $\CM$-free:
Let $s = p = 2$ with the following Auslander-Reiten quiver:
$$
%=====================================6 a=3
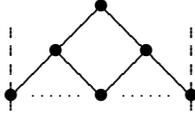
\hbox{\beginpicture
\setcoordinatesystem units <.5cm,.5cm>
\multiput{$\bullet$} at 0 0  1 1  2 0  2 2  3 1  4 0 /
\plot 0 0  2 2  4 0 /
\plot 1 1  2 0  3 1 /
\setdots <1mm>
\plot 0.5 0  1.5 0 /
\plot 2.5 0  3.5 0 /
\setdashes <1mm>
\plot 0 -.3  0 1.5 /
\plot 4 -.3  4 1.5 /
\endpicture}
$$
(thus with Kupisch series $(3,2)$). The resolution quiver looks as follows:
$$
%====================================
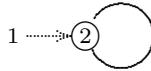
\hbox{\beginpicture
\setcoordinatesystem units <.8cm,.8cm>
\multiput{} at 1.8 .5  1.8 -.5 /
\put{$\ssize 1$} at 0 0
\put{$\ssize 2$} at 1 0
\put{$\bigcirc$} at 1 0
\circulararc 300 degrees from 1.1 -.2  center at 1.5 0 
\setdots <.4mm>
\arr{0.2 0}{0.8 0}
\endpicture} 
$$
thus there is a black loop (thus a black cycle). However, $\Lambda$
is $\CM$-free, since it has finite global dimension (the global dimension
is $2$).
   \medskip
{\bf Example 2.} A Nakayama algebra with a red loop and no other cycles 
in the resolution quiver 
which has finite global dimension: Consider the Kupisch series $(4,3,2)$:
$$
%=====================================6 a=3
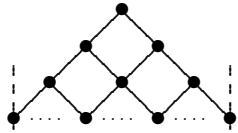
\hbox{\beginpicture
\setcoordinatesystem units <.4cm,.4cm>
\multiput{$\bullet$} at 0 0  1 1  2 0  2 2  3 1  4 0  3 3  4 2  5 1  6 0 /
\plot 0 0  3 3  6 0 /
\plot 1 1  2 0  3 1  4 0  5 1 /
\plot 2 2  3 1  4 2 /
\setdots <1mm>
\plot 0.5 0  1.5 0 /
\plot 2.5 0  3.5 0 /
\plot 4.5 0  5.5 0 /
\setdashes <1mm>
\plot 0 -.3  0 1.5 /
\plot 6 -.3  6 1.5 /
\endpicture}
$$
we have $p = 2$ and $s=3$.
The resolution quiver looks as follows:
$$
%====================================
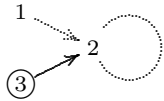
\hbox{\beginpicture
\setcoordinatesystem units <.8cm,.8cm>
\put{} at 1.8 0
\put{$\ssize 3$} at 0 -.5
\put{$\ssize 1$} at 0 .5
\put{$\ssize 2$} at 1 0
\put{$\bigcirc$} at 0 -.5
\arr{0.2 -.4}{0.8 -.1}
\setdots <.4mm>
\arr{0.2 0.4}{0.8 0.1}
\circulararc 300 degrees from 1.1 -.2  center at 1.5 0 
\endpicture} 
$$
now we deal with a red loop (thus with a cycle which is not black). As it
should be, $\Lambda$ is $\CM$-free, but it is even of finite global dimension
(the global dimension is again $2$), thus all the modules have finite
projective dimension.
	   \medskip 
{\bf Example 3.} An algebra with black and red loops in the resolution quiver
which is $\CM$-free: take the  Kupisch series $(3,3,2)$:
$$
%=====================================6 a=3
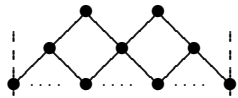
\hbox{\beginpicture
\setcoordinatesystem units <.4cm,.4cm>
\multiput{$\bullet$} at 0 0  1 1  2 0  2 2  3 1  4 0   4 2  5 1  6 0 /
\plot 0 0  2 2 /
\plot 4 2  6 0 /
\plot 1 1  2 0  3 1  4 0  5 1 /
\plot 2 2  3 1  4 2 /
\setdots <1mm>
\plot 0.5 0  1.5 0 /
\plot 2.5 0  3.5 0 /
\plot 4.5 0  5.5 0 /
\setdashes <1mm>
\plot 0 -.3  0 1.5 /
\plot 6 -.3  6 1.5 /
\endpicture}
$$
again, $p = 2$ and $s=3$.
The resolution quiver looks as follows:
$$
%====================================
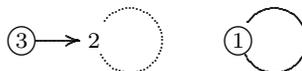
\hbox{\beginpicture
\setcoordinatesystem units <.8cm,.8cm>
\multiput{} at 0 .5  0 -.5 /
\put{$\ssize 3$} at 0 0
\put{$\ssize 2$} at 1 0
\put{$\ssize 1$} at 3 0
\multiput{$\bigcirc$} at 0 0  3 0 /
\arr{0.2 0}{0.8 0}
\circulararc 300 degrees from 3.1 -.2  center at 3.5 0 
\setdots <.4mm>
\circulararc 300 degrees from 1.1 -.2  center at 1.5 0 
\endpicture} 
$$
We have both a red loop and a black loop. The algebra $\Lambda$ has infinite
global dimension but is $\CM$-free.
       \medskip
These examples show that in Proposition 5 the assumption $p > s$ is
necessary: Namely, Example 2 exhibits an algebra $\Lambda$ of finite global dimension, thus a Gorenstein algebra, such that $R(\Lambda)$ has a red loop.
Examples 1 and 3 are $\CM$-free algebras with black cycles. Whereas
Examples 1 and 2 are algebras of finite global dimension, we do not know any
example $\Lambda$ of finite global dimension such that $R(\Lambda)$ has both
black cycles and cycles which are not black.

      \bigskip\bigskip 
%======================================================================
{\bf 8.~References.}
     \medskip
%\item{[ARS]} M. Auslander, I. Reiten, S. Smal\o. Representation Theory of
%   Artin Algebras.  Cambridge University Press (1997)
\item{[C]} X.-W. Chen. Algebras with radical square zero are either self-injective
  or $\CM$-free. Proc. Amer. Math. Soc. 140 (2012), 93-98. 
\item{[CY]} X.-W. Chen, Yu Ye. Retractions and Gorenstein homological
 properties. \newline arXiv:1206.4415. To appear in Algebras and Representation Theory.  
\item{[EJ]} E. E. Enochs, O. M. G. Jenda. Relative homological algebra.
  De Gruyter Expo\. Math\. 30   (2000).
\item{[G]} W. H. Gustafson, Global dimension in serial rings, J. Algebra 97 (1985),
  14-16.
\item {[H]} M. Hoshino: Algebras of finite self-injective dimension. 
   Proc. Amer. Math. Soc. 112 (1991), 619-622.
\item{[R1]} C.M. Ringel: 
  Representations of K-species and bimodules. J. Algebra 41 (1976), 269-302. 
\item{[R2]} C.M. Ringel: The Gorenstein projective modules for the Nakayama algebras. II.
 In preparation.
\item{[S]} D. Shen, A note on resolution quivers, in preparation.
   
%\item{[RX]} C. M. Ringel, B.-L. Xiong. On rings with radical square zero. 
%  arXiv:1112.1422

%%%%%%%%%%%%%%%%%%%%%%%%%%%%%%%%%%%%%%%%%%%%%%%%%%%%%%%%%%%%%%%%%%%%%%%%%%%%%
%%%%%%%%%%%%%%%%%%%%%%%%%%%%%%%%%%%%%%%%%%%%%%%%%%%%%%%%%%%%%%%%%%%%%%%%%%%%%
	\bigskip\bigskip
{\rmk
C. M. Ringel\par
Department of Mathematics, Shanghai Jiao Tong University \par
Shanghai 200240, P. R. China, and \par 
King Abdulaziz University, P O Box 80200\par
Jeddah, Saudi Arabia\par

%Fakult\"at f\"ur Mathematik, \par 
%Universit\"at Bielefeld, PO Box 100 131 \par 
%D-33501 Bielefeld, Germany\par 

e-mail: {\ttk ringel\@math.uni-bielefeld.de} \par
}
	\bigskip
\bye